\documentclass{amsart}
\usepackage{texdraw}

\theoremstyle{plain}
\newtheorem{thm}{Theorem}[section]

\newtheorem{lem}[thm]{Lemma}

\theoremstyle{definition}
\newtheorem{defi}[thm]{Definition}
\newtheorem{rules}[thm]{Rules}
\newtheorem{example}[thm]{Example}
\theoremstyle{remark}

\DeclareMathOperator{\wt}{wt}
\newcommand{\nc}{\newcommand}
\nc{\Uq}{U_q}
\nc{\Z}{\mathbf{Z}}
\nc{\C}{\mathbf{C}}
\nc{\Q}{\mathbf{Q}}
\nc{\op}{\oplus}
\nc{\ot}{\otimes}
\nc{\pv}{P^{\vee}}
\nc{\ali}{\alpha_i}
\nc{\B}{\mathbf{B}}
\nc{\V}{\mathbf{V}}
\nc{\La}{\Lambda}
\nc{\la}{\lambda}
\nc{\nbinom}[2]{\genfrac{}{}{0pt}{1}{{#1}}{{#2}}}
\nc{\path}{\mathcal{P}}
\nc{\fit}{\tilde{f}_i}
\nc{\eit}{\tilde{e}_i}
\nc{\Y}{\mathcal{Y}}
\nc{\A}{\mathbf{A}}
\nc{\ra}{\rightarrow}
\nc{\vep}{\varepsilon}
\nc{\vphi}{\varphi}
\nc{\h}{\mathfrak{h}}
\nc{\oP}{\bar{P}}
\nc{\tris}{
\bsegment
\move(0 0)\lvec(10 0)\lvec(10 10)\lvec(0 0)\ifill f:0.7
\esegment
}
\nc{\recs}{
\bsegment
\move(0 0)\lvec(10 0)\lvec(10 5)\lvec(0 5)\lvec(0 0)\ifill f:0.7
\esegment
}

\begin{document}

\title[Crystal graphs for $U_q(C_2^{(1)})$]
      {Crystal graphs for basic representations of
       the quantum affine algebra $U_q(C_2^{(1)})$}
\author[J. Hong and S.-J. Kang]{Jin Hong and Seok-Jin Kang$^{*}$}
\address{Department of Mathematics\\
         Seoul National University\\
         Seoul 151-742, Korea}
\thanks{$^{*}$This research was supported by KOSEF Grant
\# 98-0701-01-5-L and the Young Scientist Award, 
Korean Academy of Science and Technology}
\email{jhong@math.snu.ac.kr, sjkang@math.snu.ac.kr}

\begin{abstract}
We give a realization of crystal graphs for
basic representations of the quantum affine algebra $U_q(C_2^{(1)})$
in terms of new combinatorial objects called the Young walls.
\end{abstract}

\maketitle

\section{Introduction}

In~\cite{MR92b:17018,MR93b:17045}, Kashiwara developed the
theory of \emph{crystal bases}
for integrable modules over the quantum groups associated with symmetrizable
Kac-Moody algebras.
The crystal bases can be viewed as bases at $q=0$ and they are given a
structure of colored oriented graph, called the \emph{crystal graphs},
with arrows defined by Kashiwara operators.
The crystal graphs have many nice combinatorial features reflecting the
internal structure of integrable representations of quantum groups.
For instance, the characters of integrable representations can be
computed by counting the elements in the crystal graphs with a given weight.
Moreover, the tensor product decomposition of integrable modules into
a direct sum of irreducible submodules is equivalent to decomposing
the tensor product of crystal graphs into a disjoint union of connected
components.
Therefore, to understand the combinatorial nature of integrable
representations, it is essential to find realizations of crystal graphs
in terms of nice combinatorial objects.

In~\cite{MR91j:17021}, Misra and Miwa constructed the crystal graphs for
basic representations of quantum affine algebras $U_q(A_n^{(1)})$ using
Young diagrams with colored boxes.
Their idea was extended to construct crystal graphs for irreducible
highest weight $U_q(A_n^{(1)})$-modules of arbitrary
higher level~\cite{MR93a:17015}.
The crystal graphs constructed in \cite{MR91j:17021} 
can be parametrized by certain paths
which arise naturally in the theory of solvable lattice models.
Motivated by this observation, Kang, Kashiwara, Misra, Miwa, Nakashima
and Nakayashiki developed the theory of \emph{perfect crystals} for
general quantum affine algebras and gave a realization of crystal
graphs for irreducible highest weight modules of arbitrary higher level
in terms of \emph{paths}.
In this way, the theory of vertex models can be explained in the
language of representation theory of quantum affine algebras
and the 1-point function of the vertex model was
expressed as the quotient of the string function and the character
of the corresponding irreducible highest weight representation.

The purpose of this paper is to give a realization of crystal graphs
for basic representations of the quantum affine algebra $U_q(C_2^{(1)})$
using some new combinatorial objects which we call the \emph{Young walls}.
The Young walls consist of colored blocks that are built on
the given {\it ground-state} 
and can be viewed as
generalizations of Young diagrams.
The rules for building Young walls
are quite similar to playing with LEGO blocks and the Tetris game.
The crystal graphs for basic representations are characterized as
the set of all \emph{reduced proper Young walls}.
The weight of a Young wall can be computed easily 
by counting the number of colored blocks 
that have been added to the ground-state.
Hence the weight multiplicity is just the number of 
all reduced proper Young walls of given weight.

One can define the $U_q(C_2^{(1)})$-action on the space spanned by
all proper Young walls and can give an algorithm for finding \emph{global
basis} (or \emph{canonical basis}) associated with each reduced proper
Young wall~\cite{Ka00}.
It still remains to extend the results of this paper to quantum affine 
algebras $U_q(C_n^{(1)})$ for $n\ge 3$ and to the higher level 
integrable representations of $C_n^{(1)}$ $(n\ge 2)$. 

\vspace{5mm}
\noindent\textbf{Acknowledgments.} \ Part of this work was completed
while S.-J. Kang was visiting Yale University and 
Korea Institute for Advanced Study in the year of 1999.
He is very grateful to the faculty and staff members of Yale University 
and Korea Institute for Advanced Study for their
hospitality and support during his visit.

\vskip 5mm
\section{The quantum affine algebra of type $C_2^{(1)}$}

Let $I=\{0,1,2\}$ be the index set.
Consider the generalized Cartan matrix $A=(a_{ij})_{i,j \in I}$ 
of affine type $C_2^{(1)}$ and its Dynkin diagram:
\begin{equation*}
A = 
\begin{pmatrix}
2 & -1 & 0\\
-2 & 2 & -2\\
0 & -1 & 2
\end{pmatrix}
\quad\text{and}\quad
\text{\raisebox{-0.4\height}{
\begin{texdraw}
\drawdim mm
\fontsize{7}{7}\selectfont
\setunitscale 0.5
\textref h:C v:C
\move(0 0)
\bsegment
\move(0 0)\lcir r:2
\move(1.7 1)\lvec(13.3 1)
\move(1.7 -1)\lvec(13.3 -1)
\move(15 0)\lcir r:2
\move(16.7 1)\lvec(28.3 1)
\move(16.7 -1)\lvec(28.3 -1)
\move(30 0)\lcir r:2
\move(13 0)\lvec(10.3 2.7)\move(13 0)\lvec(10.3 -2.7)
\move(17 0)\lvec(19.7 2.7)\move(17 0)\lvec(19.7 -2.7)
\esegment
\htext(0 -5){$0$}
\htext(15 -5){$1$}
\htext(30 -5){$2$}
\end{texdraw}
}.}
\end{equation*}

Let $\pv = \Z h_0 \op \Z h_1 \op \Z h_2 \op \Z d$ be a free abelian group,
called the {\it dual weight lattice} and set 
${\mathfrak h}= \C \otimes_{\Z} \pv$.
We define the linear functionals $\alpha_i$ and $\Lambda_i$ $(i\in I)$
on ${\mathfrak h}$ by
\begin{equation*}
\begin{aligned}
\alpha_i (h_j) & = a_{ji}, \quad \alpha_i (d) = \delta_{0, i}, \\
\Lambda_i(h_j) & = \delta_{ij}, \quad \Lambda_i(d)=0 \qquad (i,j \in I).
\end{aligned}
\end{equation*}
The $\alpha_i$ (resp. $h_i$) are called the {\it simple roots} 
(resp. {\it simple coroots}) and the $\Lambda_i$ are called the
{\it fundamental weights}. 
We denote by $\Pi=\{\alpha_i | i\in I\}$ 
(resp. $\Pi^{\vee} = \{ h_i | i\in I\}$)
the set of simple roots (resp. simple coroots).

Let $c=h_0 + h_1 + h_2$ and 
$\delta = \alpha_0 + 2\alpha_1 + \alpha_2$.
Then we have 
$\alpha_i(c)=0$, $\delta(h_i)=0$ for all $i\in I$ and $\delta(d)=1$.
We call $c$ (resp. $\delta$) the {\it canonical central element}
(resp. {\it null root}).
The free abelian group $P=\Z \Lambda_0 \op \Z \Lambda_1 
\op \Z \Lambda_2 \op \Z \delta$ is called the {\it weight lattice}
and the elements of $P$ are called the {\it affine weights}. 

We denote by $q^h$ $(h\in P^{\vee})$ the basis elements
of the group algebra $\C(q)[P^{\vee}]$ with the multiplication
$q^h q^{h'} = q^{h+h'}$ $(h,h'\in \pv)$.
Set $q_0=q_2=q^2$, $q_1=q$ and $K_0=q^{2h_0}$, $K_1 = q^{h_1}$,
$K_2 = q^{2h_2}$. 
We will also use the following notations.
\begin{equation*}
[k]_i = \frac{q_i^k - q_i^{-k}}{q_i - q_i^{-1}},
\quad
[n]_i ! = \prod_{k=1}^{n} [k]_i,
\quad
\text{and}
\quad
e_i^{(n)} = e_i^n/[n]_i!,
f_i^{(n)} = f_i^n/[n]_i!.
\end{equation*}

\begin{defi}
The \emph{quantum affine algebra $U_q(C_2^{(1)})$} of type $C_2^{(1)}$
is the associative algebra with 1 over $\C(q)$ 
generated by the symbols $e_i$, $f_i$ $(i\in I)$ and $q^h$
$(h\in \pv)$ subject to the following defining relations:
\end{defi}
\begin{equation*}
\begin{aligned}
\ & q^0 = 1, \ \ q^h q^{h'} = q^{h+ h'} \quad (h, h'\in P^{\vee}),\\
\ & q^h e_i q^{-h} = q^{\ali(h)} e_i, \quad 
q^h f_i q^{-h} = q^{-\ali(h)} f_i \quad (h\in P^{\vee}, i\in I), \\
\ & e_i f_j - f_j e_i = \delta_{i,j} \frac{K_i - K_i^{-1}}{q_i - q_i^{-1}}
\quad (i,j \in I), \\
\ & e_0^2e_1 - (q^2 + q^{-2})e_0 e_1 e_0 + e_1e_0^2 = 0, \\
\ & f_0^2f_1 - (q^2 + q^{-2})f_0 f_1 f_0 + f_1f_0^2 = 0, \\
\ & e_1^3e_0 - (q^2 + 1 + q^{-2})e_1^2e_0e_1
      + (q^2 + 1 + q^{-2})e_1e_0e_1^2 - e_0e_1^3 = 0, \\
\ & f_1^3f_0 - (q^2 + 1 + q^{-2})f_1^2f_0f_1
      + (q^2 + 1 + q^{-2})f_1f_0f_1^2 - f_0f_1^3 = 0, \\
\ & e_1^3e_2 - (q^2 + 1 + q^{-2})e_1^2e_2e_1
      + (q^2 + 1 + q^{-2})e_1e_2e_1^2 - e_2e_1^3 = 0, \\
\ & f_1^3f_2 - (q^2 + 1 + q^{-2})f_1^2f_2f_1
      + (q^2 + 1 + q^{-2})f_1f_2f_1^2 - f_2f_1^3 = 0, \\
\ & e_2^2e_1 - (q^2 + q^{-2})e_2 e_1 e_2 + e_1e_2^2 = 0, \\
\ & f_2^2f_1 - (q^2 + q^{-2})f_2 f_1 f_2 + f_1f_2^2 = 0, \\
\ & e_0e_2 = e_2e_0, \quad f_0f_2 = f_2f_0.
\end{aligned}
\end{equation*}
We call $(A, \Pi, \Pi^{\vee}, P, \pv)$ the {\it Cartan datum}
associated with the quantum affine algebra $U_q(C_2^{(1)})$.

The subalgebra of $U_q(C_2^{(1)})$ generated by 
$e_i$, $f_i$, $K_i^{\pm1}$ $(i\in I)$ is denoted by 
$U'_q(C_2^{(1)})$.
It is also called the {\it quantum affine algebra of type
$C_2^{(1)}$.}
Let $\oP^{\vee} = \Z h_0 \op \Z h_1 \op \Z h_2$
and $\overline {\mathfrak h} = \C \otimes_{\Z} \oP^{\vee}$.
Consider $\alpha_i$ and $\Lambda_i$ $(i\in I)$ as linear functionals
on $\overline{\mathfrak h}$ and set 
$\oP=\Z \Lambda_0 \op \Z \Lambda_1 \op \Z \Lambda_2$.
The elements of $\oP$ are called the {\it classical weights}. 
We call $(A, \Pi, \Pi^{\vee}, \oP, \oP^{\vee})$
the {\it Cartan datum} for the quantum affine algebra
$U'_{q}(C_2^{(1)})$. 

\vskip 5mm
\section{Crystal bases}

In this section, we review the crystal basis theory for quantum affine 
algebras $U_q(C_2^{(1)})$ and $U'_q(C_2^{(1)})$.
A $U_q(C_2^{(1)})$-module (resp. 
$U'_q(C_2^{(1)})$-module) $M$ is called
{\it integrable} if 

(i) $M=\bigoplus_{\lambda \in P} M_{\lambda}$ 
(resp. $M=\bigoplus_{\lambda \in \oP} M_{\lambda}$),
where 
\begin{equation*}
M_{\la} = \{ v\in M \mid q^h v = q^{\la(h)} v \text{ for all }h\in \pv
\ (\text {resp. } h\in \oP^{\vee}) \},
\end{equation*}

(ii) for each $i\in I$, $M$ is a direct sum of 
finite dimensional irreducible $U_i$-modules, where $U_i$ denotes
the subalgebra generated by $e_i$, $f_i$, $K_i^{\pm1}$ which is
isomorphic to $\Uq(\mathfrak{sl}_2)$. 

Fix $i\in I$.
By the representation theory of $\Uq(\mathfrak{sl}_2)$, any
element $v\in M_\la$ may be written uniquely as
\begin{equation*}
v = \sum_{k\geq0} f_i^{(k)} v_k,
\end{equation*}
where $v_k \in \ker e_i \cap M_{\la+k\ali}$.
We define the endomorphisms $\eit$ and $\fit$ on $M$,
called the {\it Kashiwara operators},  by
\begin{equation*}
\eit v = \sum_{k\geq1} f_i^{(k-1)} v_k, \qquad 
\fit v = \sum_{k\geq0} f_i^{(k+1)} v_k.
\end{equation*}

Let $\A$ be the subring of $\C(q)$ consisting of the rational functions
in $q$ that are regular at $q=0$.

\begin{defi}
\hfill\\
(a) A free $\A$-submodule $L$ of an integrable $\Uq$-module $M$,
stable under $\eit$ and $\fit$, is called a \emph{crystal lattice}
if $M\cong \C(q)\ot_{\A}L$ and $L = \op_{\la\in P} L_\la$, where
$L_\la = L\cap M_\la$.\\
(b) A \emph{crystal basis} of an integrable module $M$ is a pair
$(L,B)$ such that

\hskip 3mm (i)  $L$ is a crystal lattice of $M$,

\hskip 3mm (ii) $B$ is a $\C$-basis of $L/qL$,

\hskip 3mm (iii) $B = \cup_{\la\in P} B_\la$, 
where $B_\la = B\cap(L_\la/qL_\la)$,

\hskip 3mm (iv) $\eit B\subset B\cup \{0\}$, $\fit B\subset B\cup \{0\}$,

\hskip 3mm (v) for $b, b'\in B$, $b' = \fit b$ if and only if $b = \eit b'$.

\end{defi}

\noindent
The set $B$ is given a colored oriented graph structure by defining
\raisebox{-0.1\height}{
\begin{texdraw}
\drawdim mm
\textref h:C v:C
\arrowheadsize l:2.4 w:1.1 \arrowheadtype t:F
\setunitscale 1
\htext(0 0){$b$}
\move(2 -0.5)\ravec(6 0)
\htext(10 0){$b'$}
\htext(4.5 0.8){$_i$}
\move(-1.1 -1.6)\move(11 1.5)
\end{texdraw}
}
if and only if $b' = \fit b$.
The graph $B$ is called the {\it crystal graph} of $M$ and it reflects 
the combinatorial structure of $M$.
For instance, we have $\dim_{\C(q)} M_{\lambda} =\# B_{\lambda}$ 
for all $\lambda \in P$ (or $\lambda \in \oP$). 
By extracting properties of the crystal graphs, we define the 
notion of abstract \emph{crystals} as follows. 

\begin{defi}[\cite{MR95c:17024}]
An \emph{affine crystal} (resp. \emph{classical crystal}) 
is a set $B$ together with the maps
$\wt : B \ra P$ (resp. $\wt : B \ra \oP$),
$\vep_i : B \ra \Z\cup\{-\infty\}$,
$\vphi_i : B \ra \Z\cup\{-\infty\}$,
$\eit : B \ra B\cup\{0\}$, and
$\fit : B \ra B\cup\{0\}$,
satisfying the following conditions:

\hskip 3mm (i) $\langle\wt(b),h_i\rangle = \vphi_i(b) - \vep_i(b)$ for all $b\in B$,

\hskip 3mm (ii) $\wt(\eit b) = \wt(b) + \ali$ for $b\in B$ with $\eit b\in B$, 

\hskip 3mm (iii) $\wt(\fit b) = \wt(b) - \ali$ for $b\in B$ with $\fit b\in B$,

\hskip 3mm (iv) $b' = \fit b$ if and only if $b = \eit b'$ for $b, b' \in B$,

\hskip 3mm (v) $\eit b = \fit b = 0$ if $\vep_i(b) = - \infty$.

\end{defi}

\noindent
The crystal graph of an integrable $U_q(C_2^{(1)})$-module 
(resp. $U'_q(C_2^{(1)})$-module) is an affine crystal 
(resp. a classical crystal). 

Let $B_1$ and $B_2$ be (affine or classical) crystals. 
A \emph{morphism} $\psi:B_1 \ra B_2$ of crystals is a map
$\psi:B_1 \cup \{0\} \ra B_2 \cup \{0\}$ such that 

\hskip 3mm (i) $\psi(0)=0$,

\hskip 3mm (ii) if $b\in B_1$ and $\psi(b) \in B_2$, then 
$\wt(\psi(b))= \wt(b)$, $\vep_i(\psi(b))=\vep_i(b)$, 
$\vphi_i(\psi(b)) = \vphi_i(b)$,

\hskip 3mm (iii) if $b, b'\in B_1$, $\psi(b), \psi(b') \in B_2$ and $\fit b =b'$,
then $\fit \psi(b) = \psi(b')$. 

The \emph{tensor product} $B_1 \otimes B_2$ of $B_1$ and $B_2$ 
is the set $B_1 \times B_2$ whose crystal structure is defined by 
\begin{equation*}
\begin{aligned}
\wt(b_1\ot b_2) &= \wt(b_1) + \wt(b_2),\\
\vep_i(b_1\ot b_2)
  &= \max(\vep_i(b_1), \vep_i(b_2) - \langle\wt(b_1),h_i\rangle),\\
\vphi_i(b_1\ot b_2)
  &= \max(\vphi_i(b_2), \vep_i(b_1)+\langle\wt(b_2),h_i\rangle),\\
\eit(b_1\ot b_2)
  &=
  \begin{cases}
  \eit b_1 \ot b_2 & \text{if $\vphi_i(b_1) \geq \vep_i(b_2)$,}\\
  b_1\ot\eit b_2 & \text{if $\vphi_i(b_1) < \vep_i(b_2)$,}
  \end{cases}\\
\fit(b_1\ot b_2)  
  &= 
  \begin{cases}
  \fit b_1 \ot b_2 & \text{if $\vphi_i(b_1) > \vep_i(b_2)$,}\\
  b_1\ot\fit b_2 & \text{if $\vphi_i(b_1) \leq \vep_i(b_2)$.}
  \end{cases} 
\end{aligned}
\end{equation*}

\vskip 5mm
\section{Perfect crystals}

Let $B$ be a classical crystal.
For $b\in B$, we write
$\vep(b) = \sum_i \vep_i(b)\La_i$ and
$\vphi(b) = \sum_i \vphi_i(b)\La_i$.
Note that $\wt(b) = \vphi(b) - \vep(b)$.
Set $\oP^+ = \{\la\in \oP \mid \langle \la,h_i\rangle \geq0
\text{ for all } i\in I\}$
and $\oP^+_l = \{\la\in \oP^+\mid \langle\la, c\rangle = l \}$.

\begin{defi}[\cite{MR94j:17013}]
For $l\in \Z_{> 0}$, we say that a classical crystal $B$
is a \emph{perfect crystal} of \emph{level $l$} if

(i) $B\ot B$ is connected,

(ii)  there exists some $\la_0\in\oP$ such that 
      $\wt(B)\subset \la_0 + \sum_{i\neq0} \Z_{\leq0}\alpha_i$
      and $\#(B_{\la_0}) = 1$,


(iii) for any $b\in B$, we have $\langle \vep(b), c\rangle \geq l$,

(iv)  the maps $\vep, \vphi : B^{\text{min}} = \{ b\in B \mid
      \langle \vep(b),c\rangle = l \} \ra \oP^+_l$ are bijective.

\end{defi}
\noindent
A finite dimensional $U'_q(C_2^{(1)})$-module $V$ is called a
\emph{perfect representation} of \emph{level $l$} if it has a crystal
basis $(L,B)$ such that $B$ is a perfect crystal of level $l$.

Consider the space
\begin{equation*}
\V = \C(q) v_{1,2} \op \C(q)v_{1,\bar{2}} \op \C(q)v_{2,\bar{2}}
\op\C(q) v_{2,\bar{1}} \op \C(q)v_{\bar{2},\bar{1}},
\end{equation*}
with the action of $\Uq'(C_2^{(1)})$ defined by
{\allowdisplaybreaks
\begin{align*}
e_0 v_{i,j} &=
\begin{cases}
v_{2,\bar{1}} & \text{if $(i,j) = (1,2)$,}\\
v_{\bar{2},\bar{1}} & \text{if $(i,j) = (1,\bar{2})$,}\\
0 & \text{otherwise,}
\end{cases}\\
f_0 v_{i,j} &=
\begin{cases}
v_{1,2} & \text{if $(i,j) = (2,\bar{1})$,}\\
v_{1,\bar{2}} & \text{if $(i,j) = (\bar{2},\bar{1})$,}\\
0 & \text{otherwise,}
\end{cases}\\
K_0 v_{i,j} &=
\begin{cases}
q^2v_{i,j} & \text{if $(i,j) = (2,\bar{1})$ or $(\bar{2},\bar{1})$,}\\
q^{-2}v_{i,j} & \text{if $(i,j) = (1,2)$ or $(1,\bar{2})$,}\\
v_{i,j} & \text{otherwise,}
\end{cases}\\
e_1 v_{i,j} &=
\begin{cases}
(q+q^{-1})v_{1,\bar{2}} & \text{if $(i,j) = (2,\bar{2})$,}\\
v_{2,\bar{2}} & \text{if $(i,j) = (2,\bar{1})$,}\\
0 & \text{otherwise,}
\end{cases}\\
f_1 v_{i,j} &=
\begin{cases}
v_{2,\bar{2}} & \text{if $(i,j) = (1,\bar{2})$,}\\
(q+q^{-1})v_{2,\bar{1}} & \text{if $(i,j) = (2,\bar{2})$,}\\
0 & \text{otherwise,}
\end{cases}\\
K_1 v_{i,j} &=
\begin{cases}
q^{2}v_{1,\bar{2}} & \text{if $(i,j) = (1,\bar{2})$,}\\
q^{-2}v_{2,\bar{1}} & \text{if $(i,j) = (2,\bar{1})$,}\\
v_{i,j} & \text{otherwise,}
\end{cases}\\
e_2 v_{i,j} &=
\begin{cases}
v_{1,2} & \text{if $(i,j) = (1,\bar{2})$,}\\
v_{2,\bar{1}} & \text{if $(i,j) = (\bar{2},\bar{1})$,}\\
0 & \text{otherwise,}
\end{cases}\\
f_2 v_{i,j} &=
\begin{cases}
v_{1,\bar{2}} & \text{if $(i,j) = (1,2)$,}\\
v_{\bar{2},\bar{1}} & \text{if $(i,j) = (2,\bar{1})$,}\\
0 & \text{otherwise,}
\end{cases}\\
K_2 v_{i,j} &=
\begin{cases}
q^2v_{i,j} & \text{if $(i,j) = (1,2)$ or $(2,\bar{1})$,}\\
q^{-2}v_{i,j} & \text{if $(i,j) = (1,\bar{2})$ or $(\bar{2},\bar{1})$,}\\
v_{i,j} & \text{otherwise.}
\end{cases}
\end{align*}
} 

\begin{thm}[\cite{MR94j:17013}]
The space $\V$ is an irreducible $\Uq'(C_2^{(1)})$-module 
whose crystal graph $\B$ is perfect of level 1 as is shown below: 

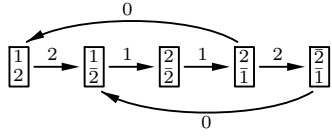
\begin{center}
\begin{texdraw}
\drawdim mm
\fontsize{7}{7}\selectfont
\textref h:C v:C
\arrowheadsize l:2.4 w:1.1 \arrowheadtype t:F
\setunitscale 1
\move(0 0)
\bsegment
\setsegscale 1.3
\move(-1 2)\lvec(1 2)\lvec(1 -2)\lvec(-1 -2)\lvec(-1 2)
\htext(0 1){$1$}\htext(0 -1){$2$}
\esegment
\move(10 0)
\bsegment
\setsegscale 1.3
\move(-1 2)\lvec(1 2)\lvec(1 -2)\lvec(-1 -2)\lvec(-1 2)
\htext(0 1){$1$}\htext(0 -1){$\bar{2}$}
\esegment
\move(20 0)
\bsegment
\setsegscale 1.3
\move(-1 2)\lvec(1 2)\lvec(1 -2)\lvec(-1 -2)\lvec(-1 2)
\htext(0 1){$2$}\htext(0 -1){$\bar{2}$}
\esegment
\move(30 0)
\bsegment
\setsegscale 1.3
\move(-1 2)\lvec(1 2)\lvec(1 -2)\lvec(-1 -2)\lvec(-1 2)
\htext(0 1){$2$}\htext(0 -1){$\bar{1}$}
\esegment
\move(40 0)
\bsegment
\setsegscale 1.3
\move(-1 2)\lvec(1 2)\lvec(1 -2)\lvec(-1 -2)\lvec(-1 2)
\htext(0 1){$\bar{2}$}\htext(0 -1){$\bar{1}$}
\esegment
\move(0 0)
\bsegment
\move(2 0)\ravec(6 0)
\move(12 0)\ravec(6 0)
\move(22 0)\ravec(6 0)
\move(32 0)\ravec(6 0)
\move(29 3.3)\clvec(23 7)(7 7)(1 3.3)
\move(3 4.3)\avec(1 3.3)
\move(39 -3.3)\clvec(33 -7)(17 -7)(11 -3.3)
\move(13 -4.3)\avec(11 -3.3)
\esegment
\move(0 0)
\bsegment
\htext(4.5 1.5){$2$}
\htext(14.5 1.5){$1$}
\htext(24.5 1.5){$1$}
\htext(34.5 1.5){$2$}
\htext(14.5 7.5){$0$}
\htext(25 -7.5){$0$}
\esegment
\move(-2 9)\move(42 -9)
\end{texdraw}
\end{center}
\end{thm}

For a dominant integral weight $\Lambda \in \oP^{+}$, we denote
by $B(\Lambda)$ the crystal graph for the irreducible highest weight 
$U'_q(C_2^{(1)})$-module $V(\Lambda)$.
As the canonical central element for $C_2^{(1)}$ is given by
$c = h_0 + h_1 + h_2$, the dominant integral weights of level one 
(i.e., those $\Lambda \in \oP^{+}$ such that $\Lambda(c)=1$)
are of the form $\Lambda = \La_i$ ($i=0,1,2$).
The level one irreducible highest weight representations are called 
the \emph{basic representations}. 

\begin{thm}[\cite{KMN1}]
There exist isomorphisms of $U'_q(C_2^{(1)})$-crystals:
\begin{equation*}
\begin{aligned}
B(\La_0) \ot \B &\cong B(\La_2), &\qquad
u_{\La_0} \ot
\text{\!\!
\raisebox{-0.4\height}{
\begin{texdraw}
\drawdim mm
\fontsize{7}{7}\selectfont
\textref h:C v:C
\setunitscale 1
\move(0 0)
\bsegment
\setsegscale 1.3
\move(-1 2)\lvec(1 2)\lvec(1 -2)\lvec(-1 -2)\lvec(-1 2)
\htext(0 1){$1$}\htext(0 -1){$2$}
\esegment
\end{texdraw}
}}
&\mapsto u_{\La_2},\\
B(\La_1) \ot \B &\cong B(\La_1), &\qquad
u_{\La_1} \ot
\text{\!\!
\raisebox{-0.4\height}{
\begin{texdraw}
\drawdim mm
\fontsize{7}{7}\selectfont
\textref h:C v:C
\setunitscale 1
\move(0 0)
\bsegment
\setsegscale 1.3
\move(-1 2)\lvec(1 2)\lvec(1 -2)\lvec(-1 -2)\lvec(-1 2)
\htext(0 1){$2$}\htext(0 -1){$\bar{2}$}
\esegment
\end{texdraw}
}}
&\mapsto u_{\La_1},\\
B(\La_2) \ot \B &\cong B(\La_0), &\qquad
u_{\La_2} \ot
\text{\!\!
\raisebox{-0.4\height}{
\begin{texdraw}
\drawdim mm
\fontsize{7}{7}\selectfont
\textref h:C v:C
\setunitscale 1
\move(0 0)
\bsegment
\setsegscale 1.3
\move(-1 2)\lvec(1 2)\lvec(1 -2)\lvec(-1 -2)\lvec(-1 2)
\htext(0 1){$\bar{2}$}\htext(0 -1){$\bar{1}$}
\esegment
\end{texdraw}
}}
&\mapsto u_{\La_0}.
\end{aligned}
\end{equation*}
\end{thm}

\vskip 5mm
\section{Path realization of crystal graphs}

Using the crystal isomorphisms given in Theorem 4.3, 
we will give a realization of the crystal graphs $B(\Lambda_i)$
in terms of paths. 
A \emph{path} is an infinite sequence 
\begin{equation*}
p=(p(k))_{k=0}^{\infty}=( \cdots, p(k+1), p(k), \cdots, p(1), p(0))
\end{equation*}
with $p(k) \in \B$ for all $k\ge 0$.
Among these, we single out the following three distinguished ones, 
called the \emph{ground-state paths}:
\begin{equation*}
\begin{aligned}
p_{\La_0} &= ( \cdots\; \nbinom{1}{2}\; \nbinom{\bar{2}}{\bar{1}}\;
                      \nbinom{1}{2}\; \nbinom{\bar{2}}{\bar{1}}\;
                      \nbinom{1}{2}\; \nbinom{\bar{2}}{\bar{1}}),\\
p_{\La_1} &= ( \cdots\; \nbinom{2}{\bar{2}}\; \nbinom{2}{\bar{2}}\;
                      \nbinom{2}{\bar{2}}\; \nbinom{2}{\bar{2}}\;
                      \nbinom{2}{\bar{2}}\; \nbinom{2}{\bar{2}}),\\
p_{\La_2} &= ( \cdots\; \nbinom{\bar{2}}{\bar{1}}\; \nbinom{1}{2}\;
                      \nbinom{\bar{2}}{\bar{1}}\; \nbinom{1}{2}\;
                      \nbinom{\bar{2}}{\bar{1}}\; \nbinom{1}{2}).
\end{aligned}
\end{equation*}
The ground-state paths are determined as the image of 
highest weight vectors in $B(\Lambda_i)$ under taking the compositions
of the inverses of the crystal isomorphisms given in Theorem 4.3. 
For example, under the isomorphism
\begin{align*}
B(\La_0)
&\cong B(\La_2)\ot\B\\
&\cong B(\La_0)\ot\B\ot\B\\
&\cong B(\La_2)\ot\B\ot\B\ot\B\\
&\cong \cdots,
\end{align*}
the highest weight vector $u_{\Lambda_0}$ is mapped onto the 
ground-state path $p_{\Lambda_0}$.

Let 
\begin{equation*}
\path(\La_i) = \{ p = (p(k))_{k=0}^{\infty} \mid
                  p(k) = p_{\La_i}(k) \text{ for all $k\gg0$} \}.
\end{equation*}
The elements of $\path(\La_i)$ are called the {\it $\Lambda_i$-paths}.
These are infinite sequence of elements from $\B$ whose tail is
identical to the ground-state path $p_{\Lambda_i}$.
We will define a crystal structure on $\path(\La_i)$ as follows. 
We need to define the action of $\fit$ and $\eit$ to each sequence.
First, replace each element of the sequence with some $0$'s and $1$'s
as shown below:
\begin{equation*}
\begin{array}{c|ccccc}
&
\text{\!\!
\raisebox{-0.4\height}{
\begin{texdraw}
\drawdim mm
\fontsize{7}{7}\selectfont
\textref h:C v:C
\setunitscale 1
\move(0 0)
\bsegment
\setsegscale 1.3
\move(-1 2)\lvec(1 2)\lvec(1 -2)\lvec(-1 -2)\lvec(-1 2)
\htext(0 1){$1$}\htext(0 -1){$2$}
\esegment
\end{texdraw}
}}
&
\text{\!\!
\raisebox{-0.4\height}{
\begin{texdraw}
\drawdim mm
\fontsize{7}{7}\selectfont
\textref h:C v:C
\setunitscale 1
\move(0 0)
\bsegment
\setsegscale 1.3
\move(-1 2)\lvec(1 2)\lvec(1 -2)\lvec(-1 -2)\lvec(-1 2)
\htext(0 1){$1$}\htext(0 -1){$\bar{2}$}
\esegment
\end{texdraw}
}}
&
\text{\!\!
\raisebox{-0.4\height}{
\begin{texdraw}
\drawdim mm
\fontsize{7}{7}\selectfont
\textref h:C v:C
\setunitscale 1
\move(0 0)
\bsegment
\setsegscale 1.3
\move(-1 2)\lvec(1 2)\lvec(1 -2)\lvec(-1 -2)\lvec(-1 2)
\htext(0 1){$2$}\htext(0 -1){$\bar{2}$}
\esegment
\end{texdraw}
}}
&
\text{\!\!
\raisebox{-0.4\height}{
\begin{texdraw}
\drawdim mm
\fontsize{7}{7}\selectfont
\textref h:C v:C
\setunitscale 1
\move(0 0)
\bsegment
\setsegscale 1.3
\move(-1 2)\lvec(1 2)\lvec(1 -2)\lvec(-1 -2)\lvec(-1 2)
\htext(0 1){$2$}\htext(0 -1){$\bar{1}$}
\esegment
\end{texdraw}
}}
&
\text{\!\!
\raisebox{-0.4\height}{
\begin{texdraw}
\drawdim mm
\fontsize{7}{7}\selectfont
\textref h:C v:C
\setunitscale 1
\move(0 0)
\bsegment
\setsegscale 1.3
\move(-1 2)\lvec(1 2)\lvec(1 -2)\lvec(-1 -2)\lvec(-1 2)
\htext(0 1){$\bar{2}$}\htext(0 -1){$\bar{1}$}
\esegment
\end{texdraw}
}}
\\[1.5mm]
\hline
\\[-3.5mm]
i=0 & 1& 1& & 0& 0\\
i=1 & &0\,0&1\,0&1\,1& \\
i=2 &0&1& &0&1
\end{array}
\end{equation*}
The blank spaces above mean to replace them with nothing.
Having done so, we remove each $0\,1$ pair occurring in the resulting
sequence of 0's and 1's.
Reading from the left, we would obtain a finite sequence of
1's followed by some finite sequence of 0's.
If we want to apply $\eit$, we apply it to the element corresponding to the
right-most 1 to obtain another sequence of elements from $\B$.
For $\fit$, we apply it to the element corresponding to the left-most 0.
Then the set $\path(\La_i)$ becomes a $U'_q(C_2^{(1)})$-crystal.
Moreover, we have:

\begin{thm}[\cite{MR94j:17013}]\label{thm:42}
For each $i=0,1,2$, there exists an isomorphism of crystals
\begin{equation*}
\path(\La_i) \cong B(\La_i).
\end{equation*}
The ground-state path $p_{\La_i}$ is mapped onto the highest weight
vector $u_{\La_i}$ under this crystal isomorphism.
\end{thm}

Hence we obtain a realization of $B(\La_i)$ in terms of paths.
Let us see some examples.

\vskip 3mm
\begin{example}\label{ex:43}
\samepage
The highest parts of the crystal graph  $B(\La_0)$ may be drawn as below.

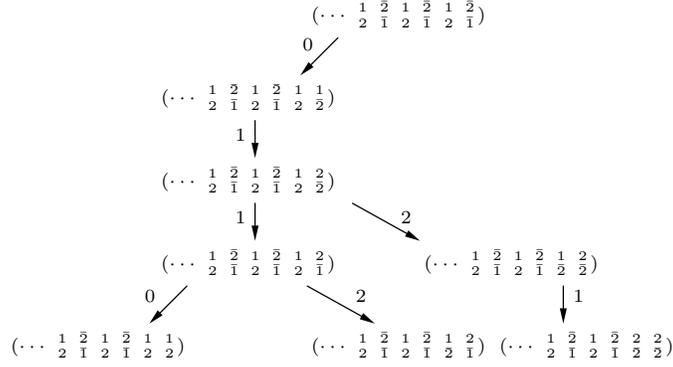
\begin{center}
\begin{texdraw}
\drawdim mm
\fontsize{7}{7}\selectfont
\textref h:C v:C
\linewd 0.2
\arrowheadsize l:2 w:0.9 \arrowheadtype t:F
\setunitscale 1
\htext(0 2){$( \cdots\; \nbinom{1}{2}\; \nbinom{\bar{2}}{\bar{1}}\;
                      \nbinom{1}{2}\; \nbinom{\bar{2}}{\bar{1}}\;
                      \nbinom{1}{2}\; \nbinom{\bar{2}}{\bar{1}})$}
\htext(-20 -9){$( \cdots\; \nbinom{1}{2}\; \nbinom{\bar{2}}{\bar{1}}\;
                      \nbinom{1}{2}\; \nbinom{\bar{2}}{\bar{1}}\;
                      \nbinom{1}{2}\; \nbinom{1}{\bar{2}})$}
\htext(-20 -20){$( \cdots\; \nbinom{1}{2}\; \nbinom{\bar{2}}{\bar{1}}\;
                      \nbinom{1}{2}\; \nbinom{\bar{2}}{\bar{1}}\;
                      \nbinom{1}{2}\; \nbinom{2}{\bar{2}})$}
\htext(-20 -31){$( \cdots\; \nbinom{1}{2}\; \nbinom{\bar{2}}{\bar{1}}\;
                      \nbinom{1}{2}\; \nbinom{\bar{2}}{\bar{1}}\;
                      \nbinom{1}{2}\; \nbinom{2}{\bar{1}})$}
\htext(15 -31){$( \cdots\; \nbinom{1}{2}\; \nbinom{\bar{2}}{\bar{1}}\;
                      \nbinom{1}{2}\; \nbinom{\bar{2}}{\bar{1}}\;
                      \nbinom{1}{\bar{2}}\; \nbinom{2}{\bar{2}})$}
\htext(-40 -42){$( \cdots\; \nbinom{1}{2}\; \nbinom{\bar{2}}{\bar{1}}\;
                      \nbinom{1}{2}\; \nbinom{\bar{2}}{\bar{1}}\;
                      \nbinom{1}{2}\; \nbinom{1}{2})$}
\htext(0 -42){$( \cdots\; \nbinom{1}{2}\; \nbinom{\bar{2}}{\bar{1}}\;
                      \nbinom{1}{2}\; \nbinom{\bar{2}}{\bar{1}}\;
                      \nbinom{1}{\bar{2}}\; \nbinom{2}{\bar{1}})$}
\htext(25 -42){$( \cdots\; \nbinom{1}{2}\; \nbinom{\bar{2}}{\bar{1}}\;
                      \nbinom{1}{2}\; \nbinom{\bar{2}}{\bar{1}}\;
                      \nbinom{2}{\bar{2}}\; \nbinom{2}{\bar{2}})$}
\move(-8 -1)\ravec(-5 -5)
\move(-19 -12)\ravec(0 -5)
\move(-19 -23)\ravec(0 -5)
\move(-28 -34)\ravec(-5 -5)
\move(-12 -34)\ravec(9 -5)
\move(-6 -23)\ravec(9 -5)
\move(22 -34)\ravec(0 -5)
\htext(-12 -2){$0$}
\htext(-21 -14){$1$}
\htext(-21 -25){$1$}
\htext(1 -25){$2$}
\htext(-33 -35.5){$0$}
\htext(-5 -35.5){$2$}
\htext(24 -35.5){$1$}
\move(37 -44.5)\move(-52 4.5)
\end{texdraw}
\end{center}
\end{example}

\vskip 5mm
\begin{example}\label{ex:32}
The highest parts of the crystal graph  $B(\La_1)$ may be drawn as below.

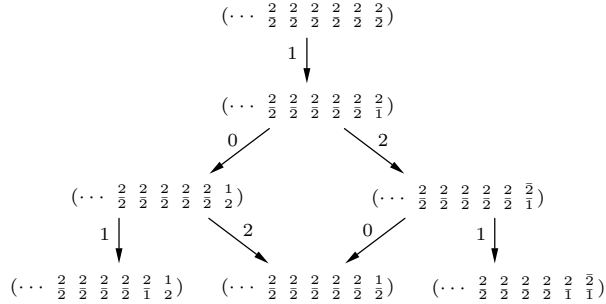
\begin{center}
\begin{texdraw}
\drawdim mm
\fontsize{7}{7}\selectfont
\textref h:C v:C
\linewd 0.2
\arrowheadsize l:2 w:0.9 \arrowheadtype t:F
\setunitscale 1
\htext(0 2){$( \cdots\; \nbinom{2}{\bar{2}}\; \nbinom{2}{\bar{2}}\;
                      \nbinom{2}{\bar{2}}\; \nbinom{2}{\bar{2}}\;
                      \nbinom{2}{\bar{2}}\; \nbinom{2}{\bar{2}})$}
\htext(0 -10){$( \cdots\; \nbinom{2}{\bar{2}}\; \nbinom{2}{\bar{2}}\;
                      \nbinom{2}{\bar{2}}\; \nbinom{2}{\bar{2}}\;
                      \nbinom{2}{\bar{2}}\; \nbinom{2}{\bar{1}})$}
\htext(-20 -22){$( \cdots\; \nbinom{2}{\bar{2}}\; \nbinom{2}{\bar{2}}\;
                      \nbinom{2}{\bar{2}}\; \nbinom{2}{\bar{2}}\;
                      \nbinom{2}{\bar{2}}\; \nbinom{1}{2})$}
\htext(20 -22){$( \cdots\; \nbinom{2}{\bar{2}}\; \nbinom{2}{\bar{2}}\;
                      \nbinom{2}{\bar{2}}\; \nbinom{2}{\bar{2}}\;
                      \nbinom{2}{\bar{2}}\; \nbinom{\bar{2}}{\bar{1}})$}
\htext(-28 -34){$( \cdots\; \nbinom{2}{\bar{2}}\; \nbinom{2}{\bar{2}}\;
                      \nbinom{2}{\bar{2}}\; \nbinom{2}{\bar{2}}\;
                      \nbinom{2}{\bar{1}}\; \nbinom{1}{2})$}
\htext(0 -34){$( \cdots\; \nbinom{2}{\bar{2}}\; \nbinom{2}{\bar{2}}\;
                      \nbinom{2}{\bar{2}}\; \nbinom{2}{\bar{2}}\;
                      \nbinom{2}{\bar{2}}\; \nbinom{1}{\bar{2}})$}
\htext(28 -34){$( \cdots\; \nbinom{2}{\bar{2}}\; \nbinom{2}{\bar{2}}\;
                      \nbinom{2}{\bar{2}}\; \nbinom{2}{\bar{2}}\;
                      \nbinom{2}{\bar{1}}\; \nbinom{\bar{2}}{\bar{1}})$}
\move(0 -1)\ravec(0 -6)\htext(-2 -3){$1$}
\move(-5 -13)\ravec(-8 -6)\htext(-10 -14.5){$0$}
\move(5 -13)\ravec(8 -6)\htext(10 -14.5){$2$}
\move(-13 -25)\ravec(8 -6)\htext(-8 -26.5){$2$}
\move(13 -25)\ravec(-8 -6)\htext(8 -26.5){$0$}
\move(25 -25)\ravec(0 -6)\htext(23 -27){$1$}
\move(-25 -25)\ravec(0 -6)\htext(-27 -27){$1$}
\move(40 -37)\move(-40 4)
\end{texdraw}
\end{center}
\end{example}

\vskip 5mm
\section{The Young walls}
The main purpose of this paper is to give a realization of crystal graph
$B(\La_i)$ using new combinatorial objects called the {\it Young walls}.
In this section, we will explain the notion of Young walls. 

The Young walls will be built of three kinds of \emph{blocks};
the $0$-block (\,%
\raisebox{-0.3\height}{%
\begin{texdraw}
\drawdim mm
\setunitscale 0.45
\fontsize{7}{7}\selectfont
\textref h:C v:C
\move(0 0)\lvec(10 0)\lvec(10 10)\lvec(0 10)\lvec(0 0)
\move(10 0)\lvec(12.5 2.5)\lvec(12.5 12.5)\lvec(2.5 12.5)\lvec(0 10)
\move(10 10)\lvec(12.5 12.5)
\htext(5 5){$0$}
\end{texdraw}%
}\,),
the $1$-block (\,%
\raisebox{-0.3\height}{%
\begin{texdraw}
\drawdim mm
\setunitscale 0.45
\fontsize{7}{7}\selectfont
\textref h:C v:C
\move(0 0)\lvec(10 0)\lvec(10 5)\lvec(0 5)\lvec(0 0)
\move(10 0)\lvec(15 5)\lvec(15 10)\lvec(5 10)\lvec(0 5)
\move(10 5)\lvec(15 10)
\htext(5 2.5){$1$}
\end{texdraw}%
}\,),
and the $2$-block (\,%
\raisebox{-0.3\height}{%
\begin{texdraw}
\drawdim mm
\setunitscale 0.45
\fontsize{7}{7}\selectfont
\textref h:C v:C
\move(0 0)\lvec(10 0)\lvec(10 10)\lvec(0 10)\lvec(0 0)
\move(10 0)\lvec(12.5 2.5)\lvec(12.5 12.5)\lvec(2.5 12.5)\lvec(0 10)
\move(10 10)\lvec(12.5 12.5)
\htext(5 5){$2$}
\end{texdraw}%
}\,).
They are supposed to be colored by elements from the index set $I$
and we do not allow rotations of the blocks.
The $0$-block and $2$-block are of unit height, unit width, and half thickness.
The $1$-block is of half height, unit width, and unit thickness.
With these blocks, we will build a wall of unit thickness, extending
infinitely to the left, like playing with LEGO blocks.
There will be many rules we must adhere to in building the wall.
The base of the wall may not be arbitrary, but must be chosen from
one of the following three.\\[5mm]
\begin{tabular}{rcl}
$Y_{\La_0}$ & $=$ &
\raisebox{-0.4\height}{
\begin{texdraw}
\drawdim mm
\setunitscale 0.5
\fontsize{7}{7}\selectfont
\textref h:C v:C
\move(0 0)\lvec(-60 0)\move(0 10)\lvec(-60 10)
\move(-57.5 12.5)\lvec(2.5 12.5)\lvec(2.5 2.5)\lvec(0 0)
\move(0 0) \bsegment \move(0 0)\lvec(0 10)\lvec(2.5 12.5) \esegment
\move(-10 0) \bsegment \move(0 0)\lvec(0 10)\lvec(2.5 12.5) \esegment
\move(-20 0) \bsegment \move(0 0)\lvec(0 10)\lvec(2.5 12.5) \esegment
\move(-30 0) \bsegment \move(0 0)\lvec(0 10)\lvec(2.5 12.5) \esegment
\move(-40 0) \bsegment \move(0 0)\lvec(0 10)\lvec(2.5 12.5) \esegment
\move(-50 0) \bsegment \move(0 0)\lvec(0 10)\lvec(2.5 12.5) \esegment
\move(-60 0) \bsegment \move(0 0)\lvec(0 10)\lvec(2.5 12.5) \esegment
\move(0 0)
\bsegment
\lpatt(0.3 1)
\move(0 0)\lvec(-2.5 -2.5)\lvec(-67.5 -2.5)
\move(-60 0)\rlvec(-5 0)\move(-60 10)\rlvec(-5 0)\move(-57.5 12.5)\rlvec(-5 0)
\move(0 0)\rlvec(-2.5 -2.5)
\move(-10 0)\rlvec(-2.5 -2.5)
\move(-20 0)\rlvec(-2.5 -2.5)
\move(-30 0)\rlvec(-2.5 -2.5)
\move(-40 0)\rlvec(-2.5 -2.5)
\move(-50 0)\rlvec(-2.5 -2.5)
\move(-60 0)\rlvec(-2.5 -2.5)
\esegment
\move(0 0)
\bsegment
\htext(-5 5){$2$}
\htext(-15 5){$0$}
\htext(-25 5){$2$}
\htext(-35 5){$0$}
\htext(-45 5){$2$}
\htext(-55 5){$0$}
\esegment
\end{texdraw}
}
\\[4mm]
$Y_{\La_1}$ & $=$ &
\raisebox{-0.4\height}{
\begin{texdraw}
\drawdim mm
\setunitscale 0.5
\fontsize{7}{7}\selectfont
\textref h:C v:C
\move(-60 0)\lvec(0 0)\move(0 5)\lvec(-60 5)
\move(0 0)\lvec(5 5)\lvec(5 10)\lvec(-55 10)
\move(0 0)\bsegment\move(0 0)\lvec(0 5)\lvec(5 10)\esegment
\move(-10 0)\bsegment\move(0 0)\lvec(0 5)\lvec(5 10)\esegment
\move(-20 0)\bsegment\move(0 0)\lvec(0 5)\lvec(5 10)\esegment
\move(-30 0)\bsegment\move(0 0)\lvec(0 5)\lvec(5 10)\esegment
\move(-40 0)\bsegment\move(0 0)\lvec(0 5)\lvec(5 10)\esegment
\move(-50 0)\bsegment\move(0 0)\lvec(0 5)\lvec(5 10)\esegment
\move(-60 0)\bsegment\move(0 0)\lvec(0 5)\lvec(5 10)\esegment
\move(0 0)
\bsegment
\lpatt(0.3 1)
\move(-60 0)\rlvec(-5 0)\move(-60 5)\rlvec(-5 0)\move(-55 10)\rlvec(-5 0)
\esegment
\move(0 0)
\bsegment
\htext(-5 2.5){$1$}
\htext(-15 2.5){$1$}
\htext(-25 2.5){$1$}
\htext(-35 2.5){$1$}
\htext(-45 2.5){$1$}
\htext(-55 2.5){$1$}
\esegment
\end{texdraw}
}
\\[4mm]
$Y_{\La_2}$ & $=$ &
\raisebox{-0.4\height}{
\begin{texdraw}
\drawdim mm
\setunitscale 0.5
\fontsize{7}{7}\selectfont
\textref h:C v:C
\move(0 0)\lvec(-60 0)\move(0 10)\lvec(-60 10)
\move(-57.5 12.5)\lvec(2.5 12.5)\lvec(2.5 2.5)\lvec(0 0)
\move(0 0) \bsegment \move(0 0)\lvec(0 10)\lvec(2.5 12.5) \esegment
\move(-10 0) \bsegment \move(0 0)\lvec(0 10)\lvec(2.5 12.5) \esegment
\move(-20 0) \bsegment \move(0 0)\lvec(0 10)\lvec(2.5 12.5) \esegment
\move(-30 0) \bsegment \move(0 0)\lvec(0 10)\lvec(2.5 12.5) \esegment
\move(-40 0) \bsegment \move(0 0)\lvec(0 10)\lvec(2.5 12.5) \esegment
\move(-50 0) \bsegment \move(0 0)\lvec(0 10)\lvec(2.5 12.5) \esegment
\move(-60 0) \bsegment \move(0 0)\lvec(0 10)\lvec(2.5 12.5) \esegment
\move(0 0)
\bsegment
\lpatt(0.3 1)
\move(0 0)\lvec(-2.5 -2.5)\lvec(-67.5 -2.5)
\move(-60 0)\rlvec(-5 0)\move(-60 10)\rlvec(-5 0)\move(-57.5 12.5)\rlvec(-5 0)
\move(0 0)\rlvec(-2.5 -2.5)
\move(-10 0)\rlvec(-2.5 -2.5)
\move(-20 0)\rlvec(-2.5 -2.5)
\move(-30 0)\rlvec(-2.5 -2.5)
\move(-40 0)\rlvec(-2.5 -2.5)
\move(-50 0)\rlvec(-2.5 -2.5)
\move(-60 0)\rlvec(-2.5 -2.5)
\esegment
\move(0 0)
\bsegment
\htext(-5 5){$0$}
\htext(-15 5){$2$}
\htext(-25 5){$0$}
\htext(-35 5){$2$}
\htext(-45 5){$0$}
\htext(-55 5){$2$}
\esegment
\end{texdraw}
}
\end{tabular}\\[5mm]
The drawings are meant to extend infinitely to the left.
The dotted lines in the front parts of $Y_{\La_0}$ and $Y_{\La_2}$ signify
where other blocks that will build up the wall of unit thickness may be placed.
These will be called the \emph{ground-state walls}.
As it is quite awkward drawing these, and since we can't see the blocks
lying to the back, we will simplify the drawing as follows.
The drawing on the right shows an example.
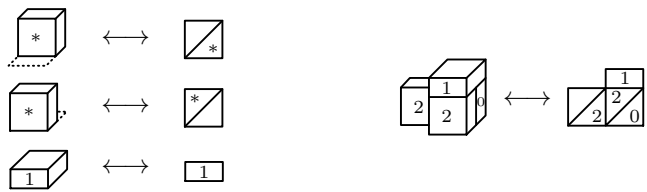
\vspace{2mm}

\begin{center}
\begin{minipage}{0.7\linewidth}
\begin{tabular}{rcl}
\raisebox{-0.4\height}{
\begin{texdraw}
\drawdim mm
\setunitscale 0.5
\fontsize{7}{7}\selectfont
\textref h:C v:C
\move(-10 0)\lvec(0 0)\lvec(0 10)\lvec(-10 10)\lvec(-10 0)
\move(0 0)\lvec(2.5 2.5)\lvec(2.5 12.5)\lvec(-7.5 12.5)\lvec(-10 10)
\move(0 10)\lvec(2.5 12.5)
\lpatt(0.3 1)
\move(0 0)\lvec(-2.5 -2.5)\lvec(-12.5 -2.5)\lvec(-10 0)
\htext(-5 5){$*$}
\end{texdraw}
}
& $\longleftrightarrow$ &
\raisebox{-0.4\height}{
\begin{texdraw}
\drawdim mm
\setunitscale 0.5
\fontsize{7}{7}\selectfont
\textref h:C v:C
\move(-10 0)\lvec(0 0)\lvec(0 10)\lvec(-10 10)\lvec(-10 0)
\move(0 10)\lvec(-10 0)
\htext(-2.5 2.5){$*$}
\end{texdraw}
}\\[4mm]
\raisebox{-0.4\height}{
\begin{texdraw}
\drawdim mm
\setunitscale 0.5
\fontsize{7}{7}\selectfont
\textref h:C v:C
\move(-10 0)\lvec(0 0)\lvec(0 10)\lvec(-10 10)\lvec(-10 0)
\move(0 0)\lvec(2.5 2.5)\lvec(2.5 12.5)\lvec(-7.5 12.5)\lvec(-10 10)
\move(0 10)\lvec(2.5 12.5)
\lpatt(0.3 1)
\move(2.5 2.5)\lvec(5 5)\lvec(2.5 5)
\htext(-5 5){$*$}
\end{texdraw}
}
& $\longleftrightarrow$ &
\raisebox{-0.4\height}{
\begin{texdraw}
\drawdim mm
\setunitscale 0.5
\fontsize{7}{7}\selectfont
\textref h:C v:C
\move(-10 0)\lvec(0 0)\lvec(0 10)\lvec(-10 10)\lvec(-10 0)
\move(0 10)\lvec(-10 0)
\htext(-7.5 7.5){$*$}
\end{texdraw}
}\\[4mm]
\raisebox{-0.4\height}{
\begin{texdraw}
\drawdim mm
\setunitscale 0.5
\fontsize{7}{7}\selectfont
\textref h:C v:C
\move(0 0)\lvec(10 0)\lvec(10 5)\lvec(0 5)\lvec(0 0)
\move(10 0)\lvec(15 5)\lvec(15 10)\lvec(5 10)\lvec(0 5)
\move(10 5)\lvec(15 10)
\htext(5 2.5){$1$}
\end{texdraw}
}
& $\longleftrightarrow$ &
\raisebox{-0.4\height}{
\begin{texdraw}
\drawdim mm
\setunitscale 0.5
\fontsize{7}{7}\selectfont
\textref h:C v:C
\move(0 0)\lvec(10 0)\lvec(10 5)\lvec(0 5)\lvec(0 0)
\htext(5 2.5){$1$}
\end{texdraw}
}
\end{tabular}
\hfill
\raisebox{-0.4\height}{
\begin{texdraw}
\drawdim mm
\setunitscale 0.5
\fontsize{7}{7}\selectfont
\textref h:C v:C
\move(0 0)\lvec(0 15)\lvec(-10 15)\lvec(-10 0)\lvec(0 0)
\move(-10 10)\lvec(0 10)\lvec(5 15)
\move(0 0)\lvec(5 5)\lvec(5 20)
\move(0 15)\lvec(5 20)\lvec(-5 20)\lvec(-10 15)
\move(2.5 2.5)\lvec(2.5 12.5)
\move(-10 2.5)\lvec(-17.5 2.5)\lvec(-17.5 12.5)\lvec(-10 12.5)
\move(-17.5 12.5)\lvec(-15 15)\lvec(-10 15)
\htext(-5 5){$2$}
\htext(-5 12.5){$1$}
\htext(-12.5 7.5){$2$}
\htext(4.2 8.7){$_0$}
\end{texdraw}
}
$\longleftrightarrow$
\raisebox{-0.37\height}{
\begin{texdraw}
\drawdim mm
\setunitscale 0.5
\fontsize{7}{7}\selectfont
\textref h:C v:C
\move(0 0)\lvec(0 10)\lvec(-20 10)\lvec(-20 0)\lvec(0 0)
\move(-20 0)\lvec(-10 10)\lvec(-10 0)\lvec(0 10)
\move(-10 10)\lvec(-10 15)\lvec(0 15)\lvec(0 10)
\htext(-12.5 2.5){$2$}
\htext(-7.5 7.5){$2$}
\htext(-2.5 2.5){$0$}
\htext(-5 12.5){$1$}
\end{texdraw}
}
\end{minipage}\\[2mm]
\end{center}

We will now list the rules for building the wall with colored
blocks.

\vskip 3mm
\begin{rules}\label{rl:51}\hfill

(i)  The wall must be built on top of one of the ground-state walls.

(ii)  The blocks should be stacked in columns.
 No block may be placed on top of a column of half thickness.

(iii) The top and the bottom of the $0$-block and $2$-block may only
      touch a $1$-block.
The sides of the $0$-block may only touch a $2$-block and vice versa.

(iv)  Placement of $0$-block and $2$-block to either the front
      or the back of the wall must be consistent for each column.
      The $0$-block and $2$-block should always be placed at heights
      which are multiples of the unit length.

(v)   Stacking more than two $1$-blocks consecutively on top of another
      is not allowed.

(vi)  Except for the right-most column, there should be no free space to
      the right of any block.

(vii) In the right-most column of a wall built on $Y_{\La_1}$,
      the $0$-block should be placed to the front and the $2$-block
      should be placed to the back.
\end{rules}

\noindent
We will give some examples illustrating the rules for building the walls.
 From now on, the ground-state wall extending infinitely to the left will be
omitted and what remains will be shaded in the drawings.

\vskip 5mm
\begin{example}\label{ex:51}
\samepage
Good walls.\\[4mm]
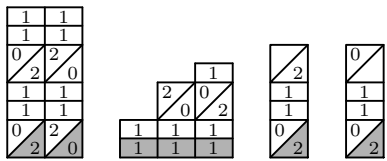
\begin{texdraw}
\drawdim mm
\setunitscale 0.5
\fontsize{7}{7}\selectfont
\textref h:C v:C
\move(0 0)\tris
\move(10 0)\tris
\move(30 0)\recs
\move(40 0)\recs
\move(50 0)\recs
\move(70 0)\tris
\move(90 0)\tris
\move(0 0)
\bsegment
\move(0 0)\rlvec(20 0)
\move(0 10)\rlvec(20 0)
\move(0 15)\rlvec(20 0)
\move(0 20)\rlvec(20 0)
\move(0 30)\rlvec(20 0)
\move(0 35)\rlvec(20 0)
\move(0 40)\rlvec(20 0)
\move(0 0)\rlvec(0 40)
\move(10 0)\rlvec(0 40)
\move(20 0)\rlvec(0 40)
\move(0 0)\rlvec(10 10)
\move(10 0)\rlvec(10 10)
\move(0 20)\rlvec(10 10)
\move(10 20)\rlvec(10 10)
\htext(2.5 7.5){$0$}
\htext(7.5 2.5){$2$}
\htext(12.5 7.5){$2$}
\htext(17.5 2.5){$0$}
\htext(2.5 27.5){$0$}
\htext(7.5 22.5){$2$}
\htext(12.5 27.5){$2$}
\htext(17.5 22.5){$0$}
\htext(5 12.5){$1$}
\htext(15 12.5){$1$}
\htext(5 17.5){$1$}
\htext(15 17.5){$1$}
\htext(5 32.5){$1$}
\htext(15 32.5){$1$}
\htext(5 37.5){$1$}
\htext(15 37.5){$1$}
\esegment
\move(60 0)
\bsegment
\move(0 0)\rlvec(-30 0)
\move(0 5)\rlvec(-30 0)
\move(0 10)\rlvec(-20 0)
\move(0 20)\rlvec(-20 0)
\move(0 25)\rlvec(-10 0)
\move(0 0)\rlvec(0 25)
\move(-10 0)\rlvec(0 25)
\move(-20 0)\rlvec(0 20)
\move(-30 0)\rlvec(0 5)
\move(-10 10)\rlvec(10 10)
\move(-20 10)\rlvec(10 10)
\move(-30 5)\lvec(-30 10)\lvec(-20 10)
\htext(-5 2.5){$1$}
\htext(-15 2.5){$1$}
\htext(-5 7.5){$1$}
\htext(-15 7.5){$1$}
\htext(-25 2.5){$1$}
\htext(-25 7.5){$1$}
\htext(-5 22.5){$1$}
\htext(-7.5 17.5){$0$}
\htext(-2.5 12.5){$2$}
\htext(-17.5 17.5){$2$}
\htext(-12.5 12.5){$0$}
\esegment
\move(80 0)
\bsegment
\move(0 0)\lvec(0 30)\lvec(-10 30)\lvec(-10 0)\lvec(0 0)
\move(-10 0)\lvec(0 10)\lvec(-10 10)
\move(0 15)\lvec(-10 15)
\move(0 20)\lvec(-10 20)\lvec(0 30)
\htext(-2.5 2.5){$2$}
\htext(-7.5 7.5){$0$}
\htext(-5 12.5){$1$}
\htext(-5 17.5){$1$}
\htext(-2.5 22.5){$2$}
\esegment
\move(100 0)
\bsegment
\move(0 0)\lvec(0 30)\lvec(-10 30)\lvec(-10 0)\lvec(0 0)
\move(-10 0)\lvec(0 10)\lvec(-10 10)
\move(0 15)\lvec(-10 15)
\move(0 20)\lvec(-10 20)\lvec(0 30)
\htext(-2.5 2.5){$2$}
\htext(-7.5 7.5){$0$}
\htext(-5 12.5){$1$}
\htext(-5 17.5){$1$}
\htext(-7.5 27.5){$0$}
\esegment
\end{texdraw}
\end{example}

\vskip 5mm

\begin{example}
\samepage
Bad walls.\\[4mm]
\begin{texdraw}
\drawdim mm
\setunitscale 0.5
\fontsize{7}{7}\selectfont
\textref h:C v:C
\move(-20 0)\tris
\move(-10 0)\tris
\move(20 0)\recs
\move(10 0)\recs
\move(40 0)\recs
\move(50 0)\recs
\move(70 0)\tris
\move(90 0)\recs
\move(110 0)\tris
\move(130 0)\recs
\move(150 0)\tris
\move(170 0)\tris
\move(180 0)\tris
\move(200 0)\recs
\move(0 0)
\bsegment
\move(-20 0)\lvec(0 0)\lvec(0 10)\lvec(-20 10)\lvec(-20 0)\lvec(-10 10)
\lvec(-10 0)\lvec(0 10)
\move(-5 10)\lvec(-5 15)\lvec(-15 15)\lvec(-15 10)
\htext(-10 12.5){$1$}
\htext(-2.5 2.5){$2$}
\htext(-7.5 7.5){$0$}
\htext(-12.5 2.5){$0$}
\htext(-17.5 7.5){$2$}
\esegment
\move(30 0)
\bsegment
\move(0 0)\lvec(-20 0)\lvec(-20 20)\lvec(0 20)\lvec(0 0)
\move(0 5)\lvec(-20 5)\move(0 10)\lvec(-20 10)
\move(0 20)\lvec(-10 10)\lvec(-10 20)\lvec(-20 10)
\move(-10 0)\lvec(-10 10)
\htext(-5 2.5){$1$}
\htext(-15 2.5){$1$}
\htext(-5 7.5){$1$}
\htext(-15 7.5){$1$}
\htext(-2.5 12.5){$2$}
\htext(-7.5 17.5){$2$}
\htext(-12.5 12.5){$0$}
\htext(-17.5 17.5){$0$}
\esegment
\move(60 0)
\bsegment
\move(0 0)\lvec(-20 0)\lvec(-20 20)\lvec(0 20)\lvec(0 0)
\move(0 5)\lvec(-20 5)\move(0 10)\lvec(-20 10)
\move(0 20)\lvec(-10 10)\lvec(-10 20)\lvec(-20 10)
\move(-10 0)\lvec(-10 10)
\htext(-5 2.5){$1$}
\htext(-15 2.5){$1$}
\htext(-5 7.5){$1$}
\htext(-15 7.5){$1$}
\htext(-2.5 12.5){$2$}
\htext(-7.5 17.5){$0$}
\htext(-12.5 12.5){$2$}
\htext(-17.5 17.5){$0$}
\esegment
\move(80 0)
\bsegment
\move(0 0)\lvec(0 30)\lvec(-10 30)\lvec(-10 0)\lvec(0 0)
\move(-10 0)\lvec(0 10)\lvec(-10 10)
\move(0 15)\lvec(-10 15)
\move(0 20)\lvec(-10 20)\lvec(0 30)
\htext(-2.5 2.5){$2$}
\htext(-7.5 7.5){$0$}
\htext(-5 12.5){$1$}
\htext(-5 17.5){$1$}
\htext(-7.5 27.5){$2$}
\htext(-2.5 22.5){$0$}
\esegment
\move(100 0)
\bsegment
\move(0 0)\lvec(0 15)\lvec(-10 15)\lvec(-10 0)\lvec(0 0)
\move(0 5)\lvec(-10 5)\lvec(0 15)
\htext(-5 2.5){$1$}
\htext(-2.5 7.5){$2$}
\htext(-7.5 12.5){$0$}
\esegment
\move(120 0)
\bsegment
\move(0 0)\lvec(0 20)\lvec(-10 20)\lvec(-10 0)\lvec(0 0)
\move(-10 0)\lvec(0 10)\lvec(-10 10)\lvec(0 20)
\htext(-2.5 2.5){$0$}
\htext(-7.5 7.5){$2$}
\htext(-2.5 12.5){$0$}
\htext(-7.5 17.5){$2$}
\esegment
\move(140 0)
\bsegment
\move(0 0)\lvec(0 15)\lvec(-10 15)\lvec(-10 0)\lvec(0 0)
\move(0 5)\lvec(-10 5)\move(0 10)\lvec(-10 10)
\htext(-5 2.5){$1$}\htext(-5 7.5){$1$}\htext(-5 12.5){$1$}
\esegment
\move(160 0)
\bsegment
\move(0 0)\lvec(0 15)\lvec(-10 15)\lvec(-10 0)\lvec(0 0)
\move(-10 10)\lvec(0 10)\lvec(-10 0)
\htext(-5 12.5){$1$}
\htext(-2.5 2.5){$2$}
\esegment
\move(190 0)
\bsegment
\move(0 0)\lvec(-20 0)\lvec(-20 10)\lvec(0 10)\lvec(0 0)
\move(-20 0)\lvec(-10 10)\lvec(-10 0)\lvec(0 10)
\htext(-17.5 7.5){$2$}
\htext(-12.5 2.5){$0$}
\htext(-2.5 2.5){$2$}
\esegment
\move(200 0)
\bsegment
\move(0 0)\lvec(10 0)\lvec(10 20)\lvec(0 20)\lvec(0 0)
\move(0 5)\lvec(10 5)
\move(10 10)\lvec(0 10)\lvec(10 20)
\htext(5 2.5){$1$}\htext(5 7.5){$1$}\htext(7.5 12.5){$0$}\htext(2.5 17.5){$2$}
\esegment
\end{texdraw}
\end{example}

A column in the wall is a \emph{full column} if the height
of the column is a multiple of the unit length and the top of the column
is of unit thickness.
The first wall in Example~\ref{ex:51} has two full columns.
The right-most column in the second wall is not full but the other two are.
Neither the third nor the forth wall contain full columns.

\begin{defi}\hfill
\samepage

(a) A wall satisfying Rules~\ref{rl:51} is called a \emph{Young wall} of
      \emph{ground state $\La_i$}, if it is built on the ground state wall
      $Y_{\La_i}$ and the height of the columns are weakly decreasing
      as we go to the left.

(b) A Young wall is said to be \emph{proper} if no two full columns of the wall
      are of the same height.
\end{defi}


\vskip 5mm
\section{The crystal structure}

In this section, we will define a crystal structure on the set of 
all proper Young walls. 
The action of Kashiwara operators will be described in a way similar to
playing the Tetris game.

\begin{defi}\hfill

(a) A block in a proper Young wall is \emph{removable}
if the wall remains a proper Young wall after removing the block.

(b) A place in a proper Young wall,
where one may add a block to obtain another proper Young wall
is called an {\it admissible slot}.

(c) A column in a proper Young wall is said to contain a
\emph{removable $\delta$}, if we may remove a $0$-block, a $2$-block,
and two $1$-blocks from the column in some order and still obtain a
proper Young wall.

(d) A proper Young wall is \emph{reduced} if none of its columns
contain a removable $\delta$.
\end{defi}

\noindent
The second and the last wall of Example~\ref{ex:51} are reduced.
The first wall is not proper and the third wall contains one
removable $\delta$.
The set of all reduced proper Young walls of ground-state $\La_i$
is denoted by $\Y(\La_i)$.


We now defined the Kashiwara operators $\eit$ and $\fit$.
A column in a proper Young wall is called 
\emph{$i$-removable} if the top block of that column is a removable 
$i$-block. 
A column is \emph{$i$-admissible} if the top of that column is an
$i$-admissible slot.
The action of Kashiwara operators is defined as follows.

\begin{enumerate}
\item Go through each column
      and write a 0 under each $i$-admissible column and a 1 under each
      $i$-removable column.
\item If we are dealing with the $i=1$ case, there could be columns from
      which two $1$-blocks may be removed.
      Place $1\,1$ under them.
      Under the columns that are $1$-removable and at the same time
      $1$-admissible, place $1\,0$.
      Write $0\,0$ under the columns that are twice $1$-admissible.
\item From the (half-)infinite list of $0$'s and $1$'s, 
      cancel out each $0\,1$ pair
      to obtain a finite sequence of $1$'s followed by some $0$'s (reading
      from left to right).
\item For $\eit$, remove the $i$-block corresponding to the
      right-most $1$ remaining.
      Set it to zero if no $1$ remains.
\item For $\fit$, add an $i$-block to the column corresponding
      to the left-most $0$ remaining.
      Set it to zero if no $0$ remains.
\end{enumerate}

\begin{lem}
The set of all reduced proper Young walls is stable under the Kashiwara
operators defined above. 
That is, if $Y$ is a reduced proper Young wall, then $\eit Y$ and $\fit Y$
are either reduced proper Young walls or 0 for all $i\in I$.
\end{lem}
\begin{proof}
By definition of removable blocks and admissible slots, 
it is obvious that the Kashiwara operators preserve proper Young walls. 
We need to verify whether they preserve \emph{reduced} proper Young walls.
We can think of two cases:

(i) adding a block may create a removable $\delta$, 

(ii) removing a block may make an existing $\delta$ removable.

Suppose we had added a $0$-block to create a removable $\delta$.
The rules for the Kashiwara operators show 
that the column standing to the left of the
created $\delta$ cannot be $0$-admissible.
So the following are essentially all possible ways the top parts of
the two columns may stand after the addition of the $0$-block.

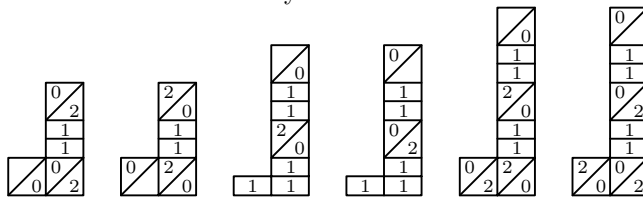
\hspace{2mm}
\begin{texdraw}
\drawdim mm
\setunitscale 0.5
\fontsize{7}{7}\selectfont
\textref h:C v:C
\move(0 0)
\bsegment
\move(10 10)\lvec(-10 10)\lvec(-10 0)\lvec(10 0)\lvec(10 30)\lvec(0 30)
\lvec(0 0)\lvec(10 10)
\move(-10 0)\lvec(0 10)
\move(0 15)\lvec(10 15)
\move(10 20)\lvec(0 20)\lvec(10 30)
\htext(-2.5 2.5){$0$} \htext(2.5 7.5){$0$} \htext(7.5 2.5){$2$}
\htext(5 12.5){$1$} \htext(5 17.5){$1$} \htext(2.5 27.5){$0$}
\htext(7.5 22.5){$2$}
\esegment
\move(30 0)
\bsegment
\move(10 10)\lvec(-10 10)\lvec(-10 0)\lvec(10 0)\lvec(10 30)\lvec(0 30)
\lvec(0 0)\lvec(10 10)
\move(-10 0)\lvec(0 10)
\move(0 15)\lvec(10 15)
\move(10 20)\lvec(0 20)\lvec(10 30)
\htext(-7.5 7.5){$0$} \htext(2.5 7.5){$2$} \htext(7.5 2.5){$0$}
\htext(5 12.5){$1$} \htext(5 17.5){$1$} \htext(2.5 27.5){$2$}
\htext(7.5 22.5){$0$}
\esegment
\move(60 0)
\bsegment
\move(10 5)\lvec(-10 5)\lvec(-10 0)\lvec(10 0)\lvec(10 40)\lvec(0 40)
\lvec(0 0)
\move(10 10)\lvec(0 10)\lvec(10 20)\lvec(0 20)
\move(10 30)\lvec(0 30)\lvec(10 40)
\move(10 25)\lvec(0 25)
\htext(-5 2.5){$1$}\htext(5 2.5){$1$}\htext(5 7.5){$1$}
\htext(2.5 17.5){$2$}\htext(7.5 12.5){$0$}
\htext(5 22.5){$1$}\htext(5 27.5){$1$}
\htext(7.5 32.5){$0$}
\esegment
\move(90 0)
\bsegment
\move(10 5)\lvec(-10 5)\lvec(-10 0)\lvec(10 0)\lvec(10 40)\lvec(0 40)
\lvec(0 0)
\move(10 10)\lvec(0 10)\lvec(10 20)\lvec(0 20)
\move(10 30)\lvec(0 30)\lvec(10 40)
\move(10 25)\lvec(0 25)
\htext(-5 2.5){$1$}\htext(5 2.5){$1$}\htext(5 7.5){$1$}
\htext(2.5 17.5){$0$}\htext(7.5 12.5){$2$}
\htext(5 22.5){$1$}\htext(5 27.5){$1$}
\htext(2.5 37.5){$0$}
\esegment
\move(120 0)
\bsegment
\move(0 0)\lvec(10 10)\lvec(-10 10)\lvec(-10 0)\lvec(10 0)\lvec(10 50)
\lvec(0 50)\lvec(0 0)
\move(-10 0)\lvec(0 10)
\move(0 15)\lvec(10 15)
\move(10 20)\lvec(0 20)\lvec(10 30)\lvec(0 30)
\move(0 35)\lvec(10 35)
\move(10 40)\lvec(0 40)\lvec(10 50)
\htext(2.5 7.5){$2$}\htext(7.5 2.5){$0$}
\htext(2.5 27.5){$2$}\htext(7.5 22.5){$0$}
\htext(7.5 42.5){$0$}
\htext(5 12.5){$1$}
\htext(5 17.5){$1$}
\htext(5 32.5){$1$}
\htext(5 37.5){$1$}
\htext(-7.5 7.5){$0$}
\htext(-2.5 2.5){$2$}
\esegment
\move(150 0)
\bsegment
\move(0 0)\lvec(10 10)\lvec(-10 10)\lvec(-10 0)\lvec(10 0)\lvec(10 50)
\lvec(0 50)\lvec(0 0)
\move(-10 0)\lvec(0 10)
\move(0 15)\lvec(10 15)
\move(10 20)\lvec(0 20)\lvec(10 30)\lvec(0 30)
\move(0 35)\lvec(10 35)
\move(10 40)\lvec(0 40)\lvec(10 50)
\htext(2.5 7.5){$0$}\htext(7.5 2.5){$2$}
\htext(2.5 27.5){$0$}\htext(7.5 22.5){$2$}
\htext(2.5 47.5){$0$}
\htext(5 12.5){$1$}
\htext(5 17.5){$1$}
\htext(5 32.5){$1$}
\htext(5 37.5){$1$}
\htext(-7.5 7.5){$2$}
\htext(-2.5 2.5){$0$}
\esegment
\end{texdraw}

\noindent
In all cases, the column containing $\delta$ already contained a removable
$\delta$ before the addition of the $0$-block, contrary to our assumption.
We may proceed similarly with other types of blocks.

Similarly, consider the case of removing an $i$-block to make an
existing $\delta$ removable.
The column standing to the right of the removed block must contain
a removable $\delta$.
The rules for the Kashiwara operators 
show that the $\delta$ column cannot be $i$-removable.
Again, we consider all possible cases subject to these two conditions
only to find that all were $\delta$-removable before the block was taken away.
\end{proof}

Therefore we have defined a crystal structure on the set $\Y(\La_i)$ of 
all reduced proper Young walls.
We now state the main result of this paper: the realization of
the crystal graph $B(\Lambda_i)$ in terms of reduced proper Young walls.

\begin{thm}
For each $i=0,1,2$, we have the isomorphism of crystals
\begin{equation*}
\Y(\La_i) \cong B(\La_i).
\end{equation*}
\end{thm}
\begin{proof}
We will prove this by giving a crystal isomorphism between $\Y(\La_i)$ and
$\path(\La_i)$.
Theorem~\ref{thm:42} tells us that this is enough.
To each element of $\Y(\La_i)$, we map an element of $\path(\La_i)$ by
reading off just the top unit cube of each column and sending
\begin{center}
\begin{tabular}{ccc}
\raisebox{-0.3\height}{
\begin{texdraw}
\drawdim mm
\setunitscale 0.5
\fontsize{7}{7}\selectfont
\textref h:C v:C
\move(0 0)\lvec(10 0)\lvec(10 10)\lvec(0 10)\lvec(0 0)\lvec(10 10)
\htext(7.5 2.5){$0$}
\end{texdraw}
}
or
\raisebox{-0.3\height}{
\begin{texdraw}
\drawdim mm
\setunitscale 0.5
\fontsize{7}{7}\selectfont
\textref h:C v:C
\move(0 0)\lvec(10 0)\lvec(10 10)\lvec(0 10)\lvec(0 0)\lvec(10 10)
\htext(2.5 7.5){$0$}
\end{texdraw}
}
&
$\longmapsto$
&
\raisebox{-0.3\height}{
\begin{texdraw}
\drawdim mm
\fontsize{7}{7}\selectfont
\textref h:C v:C
\arrowheadsize l:2.4 w:1.1 \arrowheadtype t:F
\setunitscale 1.3
\move(-1 2)\lvec(1 2)\lvec(1 -2)\lvec(-1 -2)\lvec(-1 2)
\htext(0 1){$1$}\htext(0 -1){$2$}
\end{texdraw}
},\\[3mm]
\raisebox{-0.3\height}{
\begin{texdraw}
\drawdim mm
\setunitscale 0.5
\fontsize{7}{7}\selectfont
\textref h:C v:C
\move(0 0)\lvec(10 0)\lvec(10 10)\lvec(0 10)\lvec(0 0)\lvec(10 10)
\htext(2.5 7.5){$0$}
\htext(7.5 2.5){$2$}
\end{texdraw}
}
or
\raisebox{-0.3\height}{
\begin{texdraw}
\drawdim mm
\setunitscale 0.5
\fontsize{7}{7}\selectfont
\textref h:C v:C
\move(0 0)\lvec(10 0)\lvec(10 10)\lvec(0 10)\lvec(0 0)\lvec(10 10)
\htext(7.5 2.5){$0$}
\htext(2.5 7.5){$2$}
\end{texdraw}
}
&
$\longmapsto$
&
\raisebox{-0.3\height}{
\begin{texdraw}
\drawdim mm
\fontsize{7}{7}\selectfont
\textref h:C v:C
\arrowheadsize l:2.4 w:1.1 \arrowheadtype t:F
\setunitscale 1.3
\move(-1 2)\lvec(1 2)\lvec(1 -2)\lvec(-1 -2)\lvec(-1 2)
\htext(0 1){$1$}\htext(0 -1){$\bar{2}$}
\end{texdraw}
},\\[3mm]
\raisebox{-0.3\height}{
\begin{texdraw}
\drawdim mm
\setunitscale 0.5
\fontsize{7}{7}\selectfont
\textref h:C v:C
\move(0 0)\lvec(10 0)\lvec(10 5)\lvec(0 5)\lvec(0 0)
\htext(5 2.5){$1$}
\move(0 10)
\end{texdraw}
}
&
$\longmapsto$
&
\raisebox{-0.3\height}{
\begin{texdraw}
\drawdim mm
\fontsize{7}{7}\selectfont
\textref h:C v:C
\arrowheadsize l:2.4 w:1.1 \arrowheadtype t:F
\setunitscale 1.3
\move(-1 2)\lvec(1 2)\lvec(1 -2)\lvec(-1 -2)\lvec(-1 2)
\htext(0 1){$2$}\htext(0 -1){$\bar{2}$}
\end{texdraw}
},\\[3mm]
\raisebox{-0.3\height}{
\begin{texdraw}
\drawdim mm
\setunitscale 0.5
\fontsize{7}{7}\selectfont
\textref h:C v:C
\move(0 0)\lvec(10 0)\lvec(10 10)\lvec(0 10)\lvec(0 0)
\move(0 5)\lvec(10 5)
\htext(5 2.5){$1$}
\htext(5 7.5){$1$}
\end{texdraw}
}
&
$\longmapsto$
&
\raisebox{-0.3\height}{
\begin{texdraw}
\drawdim mm
\fontsize{7}{7}\selectfont
\textref h:C v:C
\arrowheadsize l:2.4 w:1.1 \arrowheadtype t:F
\setunitscale 1.3
\move(-1 2)\lvec(1 2)\lvec(1 -2)\lvec(-1 -2)\lvec(-1 2)
\htext(0 1){$2$}\htext(0 -1){$\bar{1}$}
\end{texdraw}
},\\[3mm]
\raisebox{-0.3\height}{
\begin{texdraw}
\drawdim mm
\setunitscale 0.5
\fontsize{7}{7}\selectfont
\textref h:C v:C
\move(0 0)\lvec(10 0)\lvec(10 10)\lvec(0 10)\lvec(0 0)\lvec(10 10)
\htext(7.5 2.5){$2$}
\end{texdraw}
}
or
\raisebox{-0.3\height}{
\begin{texdraw}
\drawdim mm
\setunitscale 0.5
\fontsize{7}{7}\selectfont
\textref h:C v:C
\move(0 0)\lvec(10 0)\lvec(10 10)\lvec(0 10)\lvec(0 0)\lvec(10 10)
\htext(2.5 7.5){$2$}
\end{texdraw}
}
&
$\longmapsto$
&
\raisebox{-0.3\height}{
\begin{texdraw}
\drawdim mm
\fontsize{7}{7}\selectfont
\textref h:C v:C
\arrowheadsize l:2.4 w:1.1 \arrowheadtype t:F
\setunitscale 1.3
\move(-1 2)\lvec(1 2)\lvec(1 -2)\lvec(-1 -2)\lvec(-1 2)
\htext(0 1){$\bar{2}$}\htext(0 -1){$\bar{1}$}
\end{texdraw}
}.
\end{tabular}
\end{center}
It is clear that this map sends the ground-state walls to the appropriate
ground-state paths and that the image does indeed lie in the set
of $\La_i$-paths.
It is easy to see that this map is surjective.
Injectivity follows from the fact that the set 
$\Y(\La_i)$ consists of \emph{reduced} proper Young walls.
It remains to show that this map commutes with the action of
Kashiwara operators.

To do this, we first go through each possible case and check that
the rules for finding out which column to act on is the same
as the corresponding rules for the paths.
Let us try just one part of the $i=0$ case.
Consider the cube
\raisebox{-0.3\height}{
\begin{texdraw}
\drawdim mm
\setunitscale 0.5
\fontsize{7}{7}\selectfont
\textref h:C v:C
\move(0 0)\lvec(10 0)\lvec(10 10)\lvec(0 10)\lvec(0 0)\lvec(10 10)
\htext(2.5 7.5){$2$}
\end{texdraw}
}.
If it is a $0$-admissible slot, we would place a $0$ under the column,
which is exactly what we would do with the corresponding element
\raisebox{-0.3\height}{
\begin{texdraw}
\drawdim mm
\fontsize{7}{7}\selectfont
\textref h:C v:C
\arrowheadsize l:2.4 w:1.1 \arrowheadtype t:F
\setunitscale 1.3
\move(-1 2)\lvec(1 2)\lvec(1 -2)\lvec(-1 -2)\lvec(-1 2)
\htext(0 1){$\bar{2}$}\htext(0 -1){$\bar{1}$}
\end{texdraw}
}.
Suppose that it is not a $0$-admissible slot.
Then, the column to the right of the cube in consideration has to be a
full column of the same height.
The top cube will be
\raisebox{-0.3\height}{
\begin{texdraw}
\drawdim mm
\setunitscale 0.5
\fontsize{7}{7}\selectfont
\textref h:C v:C
\move(0 0)\lvec(10 0)\lvec(10 10)\lvec(0 10)\lvec(0 0)\lvec(10 10)
\htext(2.5 7.5){$0$}\htext(7.5 2.5){$2$}
\end{texdraw}
}.
This column is neither $0$-removable nor $0$-admissible.
So under the two columns, we would place nothing.
The corresponding two elements of $\B$ are
\raisebox{-0.3\height}{
\begin{texdraw}
\drawdim mm
\fontsize{7}{7}\selectfont
\textref h:C v:C
\arrowheadsize l:2.4 w:1.1 \arrowheadtype t:F
\setunitscale 1.3
\move(3 0)
\bsegment
\move(-1 2)\lvec(1 2)\lvec(1 -2)\lvec(-1 -2)\lvec(-1 2)
\htext(0 1){$1$}\htext(0 -1){$\bar{2}$}
\esegment
\move(0 0)
\bsegment
\move(-1 2)\lvec(1 2)\lvec(1 -2)\lvec(-1 -2)\lvec(-1 2)
\htext(0 1){$\bar{2}$}\htext(0 -1){$\bar{1}$}
\esegment
\end{texdraw}
},
under which we would place a 0 and a 1, respectively.
These cancel out when removing the $0\,1$ pair to give nothing.
So the two rules agree.
We could do the same work with other unit cubes.
The $i=1$ case is somewhat more tedious, but still possible.

Knowing that the rules for finding the column or element to act on are the
same, it suffices to check that the addition or removal of blocks match up
well with the action on the perfect crystal $\B$.
For example, adding a $2$-block to
\raisebox{-0.3\height}{
\begin{texdraw}
\drawdim mm
\setunitscale 0.5
\fontsize{7}{7}\selectfont
\textref h:C v:C
\move(0 0)\lvec(10 0)\lvec(10 10)\lvec(0 10)\lvec(0 0)\lvec(10 10)
\htext(7.5 2.5){$0$}
\end{texdraw}
}
makes it into a
\raisebox{-0.3\height}{
\begin{texdraw}
\drawdim mm
\setunitscale 0.5
\fontsize{7}{7}\selectfont
\textref h:C v:C
\move(0 0)\lvec(10 0)\lvec(10 10)\lvec(0 10)\lvec(0 0)\lvec(10 10)
\htext(7.5 2.5){$0$}
\htext(2.5 7.5){$2$}
\end{texdraw}
}.
This is in good correspondence with
\raisebox{-0.3\height}{
\begin{texdraw}
\drawdim mm
\fontsize{7}{7}\selectfont
\textref h:C v:C
\arrowheadsize l:2.4 w:1.1 \arrowheadtype t:F
\setunitscale 1
\move(0 0)
\bsegment
\setsegscale 1.3
\move(-1 2)\lvec(1 2)\lvec(1 -2)\lvec(-1 -2)\lvec(-1 2)
\htext(0 1){$1$}\htext(0 -1){$2$}
\esegment
\move(10 0)
\bsegment
\setsegscale 1.3
\move(-1 2)\lvec(1 2)\lvec(1 -2)\lvec(-1 -2)\lvec(-1 2)
\htext(0 1){$1$}\htext(0 -1){$\bar{2}$}
\esegment
\move(2 0)\ravec(6 0)
\htext(4.5 1.5){$2$}
\end{texdraw}
}.
This completes the proof.
\end{proof}

In the next examples, we redraw the crystal graphs $B(\Lambda_0)$ and
$B(\Lambda_1)$ given in Example~\ref{ex:43} and Example~\ref{ex:32}, 
this time, in terms of reduced proper Young walls.

\begin{example}
\samepage
Crystal graph of $B(\La_0)$.
\begin{center}
\begin{texdraw}
\drawdim mm
\fontsize{7}{7}\selectfont
\textref h:C v:C
\setunitscale 0.5
\htext(10 30){$Y_{\La_0}$}
\move(-20 0)\tris
\move(-20 0)
\bsegment
\move(0 0)\lvec(10 0)\lvec(10 10)\lvec(0 10)\lvec(0 0)\lvec(10 10)
\htext(2.5 7.5){$0$}\htext(7.5 2.5){$2$}
\esegment
\move(-20 -30)\tris
\move(-20 -30)
\bsegment
\move(0 0)\lvec(10 0)\lvec(10 15)\lvec(0 15)\lvec(0 0)\lvec(10 10)\lvec(0 10)
\htext(2.5 7.5){$0$}\htext(7.5 2.5){$2$}\htext(5 12.5){$1$}
\esegment
\move(-20 -65)\tris
\move(-20 -65)
\bsegment
\move(0 0)\lvec(10 0)\lvec(10 15)\lvec(0 15)\lvec(0 0)\lvec(10 10)\lvec(0 10)
\move(10 15)\lvec(10 20)\lvec(0 20)\lvec(0 15)
\htext(2.5 7.5){$0$}\htext(7.5 2.5){$2$}\htext(5 12.5){$1$}\htext(5 17.5){$1$}
\esegment
\move(86 -65)\tris
\move(76 -65)\tris
\move(86 -65)
\bsegment
\move(0 0)\lvec(10 0)\lvec(10 15)\lvec(0 15)\lvec(0 0)\lvec(10 10)\lvec(0 10)
\htext(2.5 7.5){$0$}\htext(7.5 2.5){$2$}\htext(5 12.5){$1$}
\move(0 0)\lvec(-10 0)\lvec(-10 10)\lvec(0 10)\lvec(-10 0)
\htext(-2.5 2.5){$0$}\htext(-7.5 7.5){$2$}
\esegment
\move(-50 -105)\tris
\move(-50 -105)
\bsegment
\move(0 0)\lvec(10 0)\lvec(10 30)\lvec(0 30)\lvec(0 0)\lvec(10 10)\lvec(0 10)
\move(0 15)\lvec(10 15)
\move(10 20)\lvec(0 20)\lvec(10 30)
\htext(7.5 2.5){$2$} \htext(2.5 7.5){$0$} \htext(5 12.5){$1$}
\htext(5 17.5){$1$} \htext(2.5 27.5){$0$}
\esegment
\move(10 -105)\tris
\move(0 -105)\tris
\move(10 -105)
\bsegment
\move(0 0)\lvec(10 0)\lvec(10 15)\lvec(0 15)\lvec(0 0)\lvec(10 10)\lvec(0 10)
\move(10 15)\lvec(10 20)\lvec(0 20)\lvec(0 15)
\htext(2.5 7.5){$0$}\htext(7.5 2.5){$2$}\htext(5 12.5){$1$}\htext(5 17.5){$1$}
\move(0 0)\lvec(-10 0)\lvec(-10 10)\lvec(0 10)\lvec(-10 0)
\htext(-7.5 7.5){$2$}\htext(-2.5 2.5){$0$}
\esegment
\move(100 -105)\tris
\move(90 -105)\tris
\move(100 -105)
\bsegment
\move(0 0)\lvec(10 0)\lvec(10 15)\lvec(0 15)\lvec(0 0)\lvec(10 10)\lvec(0 10)
\htext(2.5 7.5){$0$}\htext(7.5 2.5){$2$}\htext(5 12.5){$1$}
\move(-10 0)
\bsegment
\move(0 0)\lvec(10 0)\lvec(10 15)\lvec(0 15)\lvec(0 0)\lvec(10 10)\lvec(0 10)
\htext(2.5 7.5){$2$}\htext(7.5 2.5){$0$}\htext(5 12.5){$1$}
\esegment
\esegment
\move(-25 -145)\tris
\move(-35 -145)\tris
\move(-25 -145)
\bsegment
\move(0 0)\lvec(10 0)\lvec(10 30)\lvec(0 30)\lvec(0 0)\lvec(10 10)\lvec(0 10)
\move(0 15)\lvec(10 15)
\move(10 20)\lvec(0 20)\lvec(10 30)
\htext(7.5 2.5){$2$} \htext(2.5 7.5){$0$} \htext(5 12.5){$1$}
\htext(5 17.5){$1$} \htext(2.5 27.5){$0$}
\move(0 0)\lvec(-10 0)\lvec(-10 10)\lvec(0 10)\lvec(-10 0)
\htext(-7.5 7.5){$2$}
\htext(-2.5 2.5){$0$}
\esegment
\move(38 -145)\tris
\move(28 -145)\tris
\move(38 -145)
\bsegment
\move(0 0)\lvec(10 0)\lvec(10 15)\lvec(0 15)\lvec(0 0)\lvec(10 10)\lvec(0 10)
\move(10 15)\lvec(10 20)\lvec(0 20)\lvec(0 15)
\htext(2.5 7.5){$0$}\htext(7.5 2.5){$2$}\htext(5 12.5){$1$}\htext(5 17.5){$1$}
\move(0 0)\lvec(-10 0)\lvec(-10 10)\lvec(0 10)\lvec(-10 0)
\htext(-7.5 7.5){$2$}\htext(-2.5 2.5){$0$}
\move(10 20)\lvec(10 30)\lvec(0 30)\lvec(0 20)\lvec(10 30)
\htext(7.5 22.5){$2$}
\esegment
\move(62 -145)\tris
\move(72 -145)\tris
\move(82 -145)\tris
\move(72 -145)
\bsegment
\move(20 10)\lvec(-10 10)\lvec(-10 0)\lvec(20 0)\lvec(20 15)\lvec(0 15)
\lvec(0 0)\lvec(10 10)
\move(10 15)\lvec(10 0)\lvec(20 10)
\move(-10 0)\lvec(0 10)
\htext(2.5 7.5){$2$}\htext(7.5 2.5){$0$}
\htext(5 12.5){$1$}\htext(15 12.5){$1$}
\htext(12.5 7.5){$0$}\htext(17.5 2.5){$2$}
\htext(-7.5 7.5){$0$}\htext(-2.5 2.5){$2$}
\esegment
\move(113 -145)\tris
\move(103 -145)\tris
\move(113 -145)
\bsegment
\move(10 15)\lvec(-10 15)\lvec(-10 0)\lvec(10 0)\lvec(10 20)\lvec(0 20)
\lvec(0 0)\lvec(10 10)\lvec(-10 10)
\move(-10 0)\lvec(0 10)
\htext(2.5 7.5){$0$}\htext(7.5 2.5){$2$}
\htext(-7.5 7.5){$2$}\htext(-2.5 2.5){$0$}
\htext(-5 12.5){$1$}\htext(5 12.5){$1$}\htext(5 17.5){$1$}
\esegment
\move(0 0)
\bsegment
\linewd 0.45
\arrowheadsize l:3 w:1.5 \arrowheadtype t:F
\move(4 27)\avec(-6 13) \htext(-3 22){$0$}
\move(-15 -2)\avec(-15 -13) \htext(-18 -7){$1$}
\move(-15 -32)\avec(-15 -43) \htext(-18 -37){$1$}
\move(-3 -21)\ravec(74 -27)\htext(36 -32){$2$}
\move(-5 -62)\ravec(14 -15)\htext(4 -67){$2$}
\move(-25 -62)\ravec(-10 -14)\htext(-32 -67){$0$}
\move(3 -110)\ravec(-10 -14)\htext(-4 -115){$0$}
\move(-43 -110)\ravec(14 -15)\htext(-33 -115){$2$}
\move(17 -110)\ravec(14 -15)\htext(26 -115){$2$}
\move(93 -69)\ravec(0 -18)\htext(90 -77){$1$}
\move(107 -108)\ravec(0 -15)\htext(104 -115){$1$}
\move(94 -108)\ravec(-10 -15)\htext(87 -113){$0$}
\esegment
\move(0 33)
\end{texdraw}
\end{center}
\end{example}

\begin{example}
\samepage
Crystal graph of $B(\La_1)$.
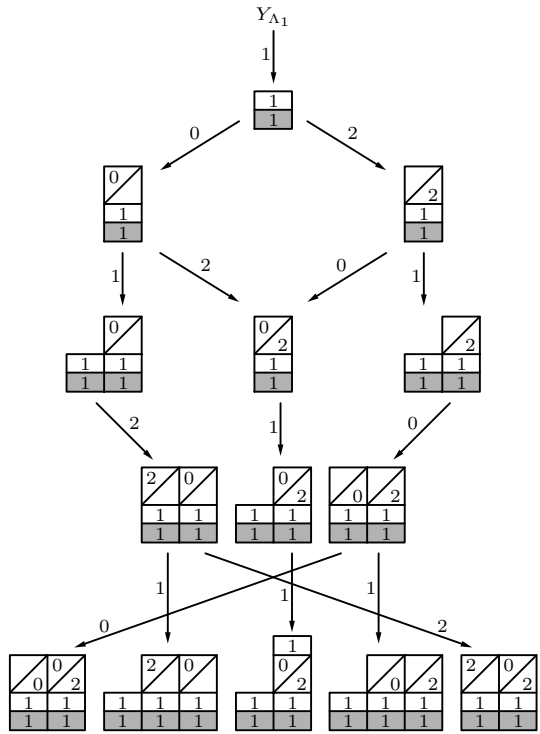
\begin{center}
\begin{texdraw}
\drawdim mm
\fontsize{7}{7}\selectfont
\textref h:C v:C
\setunitscale 0.5
\htext(5 30){$Y_{\La_1}$}
\move(0 33)
\move(0 0)\recs
\move(0 0)
\bsegment
\move(0 0)\lvec(10 0)\lvec(10 10)\lvec(0 10)\lvec(0 0)
\move(10 5)\lvec(0 5)
\htext(5 2.5){$1$}\htext(5 7.5){$1$}
\esegment
\move(-40 -30)\recs
\move(-40 -30)
\bsegment
\move(0 0)\lvec(10 0)\lvec(10 20)\lvec(0 20)\lvec(0 0)
\move(0 5)\lvec(10 5)
\move(10 10)\lvec(0 10)\lvec(10 20)
\htext(5 2.5){$1$}\htext(5 7.5){$1$}\htext(2.5 17.5){$0$}
\esegment
\move(40 -30)\recs
\move(40 -30)
\bsegment
\move(0 0)\lvec(10 0)\lvec(10 20)\lvec(0 20)\lvec(0 0)
\move(0 5)\lvec(10 5)
\move(10 10)\lvec(0 10)\lvec(10 20)
\htext(5 2.5){$1$}\htext(5 7.5){$1$}\htext(7.5 12.5){$2$}
\esegment
\move(-40 -70)\recs
\move(-50 -70)\recs
\move(-40 -70)
\bsegment
\move(10 10)\lvec(-10 10)\lvec(-10 0)\lvec(10 0)\lvec(10 20)\lvec(0 20)
\lvec(0 0)
\move(-10 5)\lvec(10 5)
\move(0 10)\lvec(10 20)
\htext(5 2.5){$1$}\htext(5 7.5){$1$}\htext(2.5 17.5){$0$}
\htext(-5 2.5){$1$}\htext(-5 7.5){$1$}
\esegment
\move(0 -70)\recs
\move(0 -70)
\bsegment
\move(0 0)\lvec(10 0)\lvec(10 20)\lvec(0 20)\lvec(0 0)
\move(0 5)\lvec(10 5)
\move(10 10)\lvec(0 10)\lvec(10 20)
\htext(5 2.5){$1$}\htext(5 7.5){$1$}\htext(7.5 12.5){$2$}\htext(2.5 17.5){$0$}
\esegment
\move(50 -70)\recs
\move(40 -70)\recs
\move(50 -70)
\bsegment
\move(10 10)\lvec(-10 10)\lvec(-10 0)\lvec(10 0)\lvec(10 20)\lvec(0 20)
\lvec(0 0)
\move(-10 5)\lvec(10 5)
\move(0 10)\lvec(10 20)
\htext(5 2.5){$1$}\htext(5 7.5){$1$}\htext(7.5 12.5){$2$}
\htext(-5 2.5){$1$}\htext(-5 7.5){$1$}
\esegment
\move(-20 -110)\recs
\move(-30 -110)\recs
\move(-20 -110)
\bsegment
\move(0 20)\lvec(0 0)\lvec(10 0)\lvec(10 20)\lvec(-10 20)\lvec(-10 0)\lvec(0 0)
\move(-10 5)\lvec(10 5)
\move(10 10)\lvec(0 10)\lvec(10 20)
\move(0 10)\lvec(-10 10)\lvec(0 20)
\htext(5 2.5){$1$}\htext(5 7.5){$1$}
\htext(-5 2.5){$1$}\htext(-5 7.5){$1$}
\htext(-7.5 17.5){$2$}\htext(2.5 17.5){$0$}
\esegment
\move(5 -110)\recs
\move(-5 -110)\recs
\move(5 -110)
\bsegment
\move(10 10)\lvec(-10 10)\lvec(-10 0)\lvec(10 0)\lvec(10 20)\lvec(0 20)
\lvec(0 0)
\move(-10 5)\lvec(10 5)
\move(0 10)\lvec(10 20)
\htext(5 2.5){$1$}\htext(5 7.5){$1$}\htext(2.5 17.5){$0$}\htext(7.5 12.5){$2$}
\htext(-5 2.5){$1$}\htext(-5 7.5){$1$}
\esegment
\move(30 -110)\recs
\move(20 -110)\recs
\move(30 -110)
\bsegment
\move(0 20)\lvec(0 0)\lvec(10 0)\lvec(10 20)\lvec(-10 20)\lvec(-10 0)\lvec(0 0)
\move(-10 5)\lvec(10 5)
\move(10 10)\lvec(0 10)\lvec(10 20)
\move(0 10)\lvec(-10 10)\lvec(0 20)
\htext(5 2.5){$1$}\htext(5 7.5){$1$}
\htext(-5 2.5){$1$}\htext(-5 7.5){$1$}
\htext(-2.5 12.5){$0$}\htext(7.5 12.5){$2$}
\esegment
\move(-55 -160)\recs
\move(-65 -160)\recs
\move(-55 -160)
\bsegment
\move(0 20)\lvec(0 0)\lvec(10 0)\lvec(10 20)\lvec(-10 20)\lvec(-10 0)\lvec(0 0)
\move(-10 5)\lvec(10 5)
\move(10 10)\lvec(0 10)\lvec(10 20)
\move(0 10)\lvec(-10 10)\lvec(0 20)
\htext(5 2.5){$1$}\htext(5 7.5){$1$}
\htext(-5 2.5){$1$}\htext(-5 7.5){$1$}
\htext(-2.5 12.5){$0$}\htext(7.5 12.5){$2$}\htext(2.5 17.5){$0$}
\esegment
\move(-20 -160)\recs
\move(-30 -160)\recs
\move(-40 -160)\recs
\move(-20 -160)
\bsegment
\move(0 20)\lvec(0 0)\lvec(10 0)\lvec(10 20)\lvec(-10 20)\lvec(-10 0)\lvec(0 0)
\move(-10 5)\lvec(10 5)
\move(10 10)\lvec(0 10)\lvec(10 20)
\move(0 10)\lvec(-10 10)\lvec(0 20)
\htext(5 2.5){$1$}\htext(5 7.5){$1$}
\htext(-5 2.5){$1$}\htext(-5 7.5){$1$}
\htext(-7.5 17.5){$2$}\htext(2.5 17.5){$0$}
\move(-10 0)\lvec(-20 0)\lvec(-20 10)\lvec(-10 10)
\move(-10 5)\lvec(-20 5)
\htext(-15 2.5){$1$}
\htext(-15 7.5){$1$}
\esegment
\move(5 -160)\recs
\move(-5 -160)\recs
\move(5 -160)
\bsegment
\move(10 10)\lvec(-10 10)\lvec(-10 0)\lvec(10 0)\lvec(10 20)\lvec(0 20)
\lvec(0 0)
\move(-10 5)\lvec(10 5)
\move(0 10)\lvec(10 20)
\htext(5 2.5){$1$}\htext(5 7.5){$1$}\htext(2.5 17.5){$0$}\htext(7.5 12.5){$2$}
\htext(-5 2.5){$1$}\htext(-5 7.5){$1$}
\move(10 20)\lvec(10 25)\lvec(0 25)\lvec(0 20)
\htext(5 22.5){$1$}
\esegment
\move(40 -160)\recs
\move(30 -160)\recs
\move(20 -160)\recs
\move(40 -160)
\bsegment
\move(0 20)\lvec(0 0)\lvec(10 0)\lvec(10 20)\lvec(-10 20)\lvec(-10 0)\lvec(0 0)
\move(-10 5)\lvec(10 5)
\move(10 10)\lvec(0 10)\lvec(10 20)
\move(0 10)\lvec(-10 10)\lvec(0 20)
\htext(5 2.5){$1$}\htext(5 7.5){$1$}
\htext(-5 2.5){$1$}\htext(-5 7.5){$1$}
\htext(-2.5 12.5){$0$}\htext(7.5 12.5){$2$}
\move(-10 0)\lvec(-20 0)\lvec(-20 10)\lvec(-10 10)
\move(-10 5)\lvec(-20 5)
\htext(-15 2.5){$1$}
\htext(-15 7.5){$1$}
\esegment
\move(65 -160)\recs
\move(55 -160)\recs
\move(65 -160)
\bsegment
\move(0 20)\lvec(0 0)\lvec(10 0)\lvec(10 20)\lvec(-10 20)\lvec(-10 0)\lvec(0 0)
\move(-10 5)\lvec(10 5)
\move(10 10)\lvec(0 10)\lvec(10 20)
\move(0 10)\lvec(-10 10)\lvec(0 20)
\htext(5 2.5){$1$}\htext(5 7.5){$1$}
\htext(-5 2.5){$1$}\htext(-5 7.5){$1$}
\htext(-7.5 17.5){$2$}\htext(2.5 17.5){$0$}
\htext(7.5 12.5){$2$}
\esegment
\move(0 0)
\bsegment
\linewd 0.45
\arrowheadsize l:3 w:1.5 \arrowheadtype t:F
\move(5 26)\ravec(0 -14)\htext(3 20){$1$}
\move(-4 2)\ravec(-21 -13)\htext(-16 -1){$0$}
\move(14 2)\ravec(21 -13)\htext(26 -1){$2$}
\move(-25 -33)\ravec(21 -13)\htext(-13 -36){$2$}
\move(35 -33)\ravec(-21 -13)\htext(23 -36){$0$}
\move(-35 -33)\ravec(0 -14)\htext(-37 -39){$1$}
\move(45 -33)\ravec(0 -14)\htext(43 -39){$1$}
\move(7 -73)\ravec(0 -14)\htext(5 -79){$1$}
\move(52 -73)\ravec(-15 -15)\htext(42 -78){$0$}
\move(-42 -73)\ravec(15 -15)\htext(-32 -78){$2$}
\move(23 -113)\ravec(-71 -25)\htext(-40 -132){$0$}
\move(-13 -113)\ravec(71 -25)\htext(50 -132){$2$}
\move(33 -113)\ravec(0 -24)\htext(31 -122){$1$}
\move(-23 -113)\ravec(0 -24)\htext(-25 -122){$1$}
\move(10 -113)\ravec(0 -21)\htext(8 -124){$1$}
\esegment
\end{texdraw}
\end{center}
\end{example}

\vskip 1cm

\begin{thebibliography}{1}

\bibitem{MR93a:17015}
M.~Jimbo, K.~C. Misra, T.~Miwa, and M.~Okado, \emph{Combinatorics of
  representations of ${U}_q(\widehat{\mathfrak{sl}}(n))$ at $q=0$}, Comm. Math.
  Phys. \textbf{136} (1991), no.~3, 543--566.

\bibitem{Ka00}
S.-J. Kang, \emph{Combinatorial representation theory of quantum affine
  algebras}, in preparation.

\bibitem{KMN1}
S.-J. Kang, M.~Kashiwara, K.~C. Misra, T.~Miwa, T.~Nakashima, and
  A.~Nakayashiki, \emph{Affine crystals and vertex models}, Int. J. Mod. Phys.
  A. \textbf{Suppl. 1A} (1992), 449--484.

\bibitem{MR94j:17013}
\bysame, \emph{Perfect crystals of quantum affine {L}ie algebras}, Duke Math.
  J. \textbf{68} (1992), no.~3, 499--607.

\bibitem{MR92b:17018}
M.~Kashiwara, \emph{Crystalizing the $q$-analogue of universal enveloping
  algebras}, Comm. Math. Phys. \textbf{133} (1990), no.~2, 249--260.

\bibitem{MR93b:17045}
\bysame, \emph{On crystal bases of the ${Q}$-analogue of universal enveloping
  algebras}, Duke Math. J. \textbf{63} (1991), no.~2, 465--516.

\bibitem{MR95c:17024}
\bysame, \emph{Crystal bases of modified quantized enveloping algebra}, Duke
  Math. J. \textbf{73} (1994), no.~2, 383--413.

\bibitem{MR91j:17021}
K.~Misra and T.~Miwa, \emph{Crystal base for the basic representation of
  ${U}_q(\widehat{\mathfrak{sl}}(n))$}, Comm. Math. Phys. \textbf{134} (1990),
  no.~1, 79--88.

\end{thebibliography}

\providecommand{\bysame}{\leavevmode\hbox to3em{\hrulefill}\thinspace}

\end{document}